\documentclass[11pt]{article}
\usepackage{cite}
\usepackage{amsmath,amssymb,amsfonts}
\usepackage{algorithm}
\usepackage[noend]{algpseudocode}
\usepackage{graphicx}
\usepackage{textcomp}
\usepackage{xcolor}
\usepackage{braket}
\usepackage{subcaption}
\usepackage{url}
\def\BibTeX{{\rm B\kern-.05em{\sc i\kern-.025em b}\kern-.08em
    T\kern-.1667em\lower.7ex\hbox{E}\kern-.125emX}}
    
\usepackage{tikz}
\usetikzlibrary{shapes.geometric, arrows}
\tikzstyle{startstop} = [rectangle, rounded corners, minimum width=3cm, minimum height=1cm,text centered, draw=black, fill=red!30]
\tikzstyle{io} = [trapezium, trapezium left angle=70, trapezium right angle=110, minimum width=3cm, minimum height=1cm, text centered, draw=black, fill=blue!30]
\tikzstyle{process} = [rectangle, minimum width=3cm, minimum height=1cm, text centered, draw=black, fill=orange!30]
\tikzstyle{decision} = [diamond, minimum width=2cm, minimum height=0.5cm, text centered, draw=black, fill=green!30]
\tikzstyle{arrow} = [thick,->,>=stealth]
\usepackage[colorlinks=true,citecolor=blue]{hyperref}
\usepackage{rotating}
\usepackage{fancyvrb}
\usepackage{changepage}
\usepackage{breakurl}

\newtheorem{definition}{Definition}

\begin{document}

\title{A Quantum Inspired Bi-level Optimization Algorithm for the First Responder Network Design Problem}
\author{
Anthony Karahalios$^{1}$\footnote{akarahal@andrew.cmu.edu}, \hspace{0.1pt}
 Sridhar Tayur$^{1}$\footnote{stayur@andrew.cmu.edu}, \hspace{0.1pt} 
 Ananth Tenneti$^{1}$\footnote{vat@andrew.cmu.edu},\\
 Amirreza Pashapour$^{2}$, \hspace{0.1pt}
 F. Sibel Salman$^{2}$, \hspace{0.1pt}
 Barış Yıldız$^{2}$ \\
$^{1}$  Carnegie Mellon University, Pittsburgh, PA 15213, USA\\
$^2$   Koç University, Istanbul, 34450, Türkiye\\
}
\maketitle

\begin{abstract}
In the aftermath of a sudden catastrophe, First Responders (FR) strive to promptly reach and rescue immobile victims. Simultaneously, other mobile individuals take roads to evacuate the affected region, access medical facilities or shelters, or reunite with their relatives. The escalated traffic congestion significantly hinders critical FR operations if they share some of the same roads. A proposal from the Turkish Ministry of Transportation and Infrastructure being discussed for implementation is to allocate a subset of road segments for use by FRs only, mark them clearly, and pre-communicate them to the citizens. For the FR paths under consideration: (i) there should exist an FR path from designated entry points to each demand point in the network, and (ii) evacuees try to leave the network (through some exit points following the selfish routing principle) in the shortest time possible when they know that certain segments are not available to them. We develop a mixed integer non-linear programming formulation for this First Responder Network Design Problem (FRNDP). We solve FRNDP using a novel hybrid quantum-classical heuristic building on the Graver Augmented Multi-Seed Algorithm (GAMA). Using the flow-balance constraints for the FR and evacuee paths, we use a Quadratic Unconstrained Binary Optimization (QUBO) model to obtain a partial Graver Bases to move between the feasible solutions of FRNDP. To efficiently explore the solution space for high-quality solutions, we develop a novel bi-level nested GAMA within GAMA: GAGA. We test GAGA on random graph instances of various sizes and instances related to an expected Istanbul earthquake. Comparing GAGA against a state-of-the-art exact algorithm for traditional formulations, we find that GAGA offers a promising alternative approach. We hope our work encourages further study of quantum (inspired) algorithms to tackle complex optimization models from other application domains.
\end{abstract}


\section{Introduction}\label{sec:intro}

Disaster Management (DM) is an important (and growing) area of study in our operations research community. DM research has been stratified into four groups: mitigation, preparedness, response, and recovery. This paper adds to the preparedness research - pre-disaster planning for a swift and effective post-disaster response - by studying the First Responder Network Design Problem (FRNDP). 
%Specifically, we aim to facilitate first-aid responders' access immediately after a disaster by identifying in the pre-disaster stage which road segments should be reserved partially or entirely for the FRs. This is accomplished by defining and solving an optimization problem that accounts for the traffic caused by the disaster-affected individuals' travels between their origin-destination (O-D) pairs.  

%On a global scale, the occurrence of natural disasters has surged tenfold since the 1960s \cite{ETR2020}. At the same time, the damage caused by these incidents has also increased. 
%As the frequency and scale of natural disasters continue to intensify with the imminent climate change, 
The complex and time-critical nature of response operations mandates meticulous pre-disaster planning. This proactive approach, employing analytical methodologies, is essential to facilitate the swift implementation of pre-identified (and pre-communicated) policies immediately when a disaster strikes. The two consecutive earthquakes of magnitudes 7.8 and 7.5 experienced in Türkiye on February 6, 2023, that affected over 15 million people in 11 provinces in southwestern Türkiye and northwestern Syria, exemplify the complexity of the required response operations. The earthquakes caused over 50,000 deaths, more than 100,000 injuries, and over 2.7 million displaced people in Türkiye alone due to 36,932 collapsed and more than 311,000 damaged buildings \cite{GRADEReport}. After the earthquakes, 271,060 personnel were deployed to the region, including 35,250 search and rescue personnel, public employees, personnel of NGOs, international search and rescue personnel, and volunteers \cite{TRReport}. Although rescue and response teams tried to reach the disaster-stricken region immediately, time lags occurred due to difficulties in accessibility, and not all cities and villages could be provided emergency services within the critical first 24 hours \cite{METUReport}. 

After a sudden-onset disaster, such as an earthquake, mobile individuals naturally attempt to evacuate the affected region, access medical facilities or shelters, and reunite with their relatives \cite{goltz2011public}. Simultaneously, first responders (FR) comprising search-and-rescue teams, firefighters, police officers, ambulance crews, debris clearance, and road restoration teams, relief aid distribution teams, and volunteers rely on the same road infrastructure to promptly reach victims in urgent need \cite{alexander2009first}. The inevitable traffic congestion caused by the mobile citizenry significantly hinders critical first-response operations. 

One of the operational planning proposals to ensure the accessibility of the FRs is to reserve one or more lanes on portions of the road network exclusively for the use of FR vehicles during the crucial initial hours following a disaster also referred to as ``golden hours" \cite{lerner2001golden}.  However, this course of action requires careful planning, as this allocation diminishes the already congested road capacity for the mobile public. Furthermore, road reservations must be decided and announced before the occurrence of a disaster. Attempting to communicate and enforce such a strategy amidst the chaos of a disaster, compounded by potential communication breakdowns, would pose an extremely challenging, if not impossible, task. Knowing which road segments should be dedicated to the FRs also has the benefit that they can be strengthened via retrofitting projects in the pre-disaster stage to withstand the disaster \cite{peeta2010, yucel2018, ZHANG2023}. The Transportation and Infrastructure Ministry of Türkiye (on August 26th, 2023) announced which highways would be reserved for FRs as part of preparation to better respond to a potential Istanbul earthquake and stated that critical structures on these roads would be strengthened. However, a detailed analysis regarding the urban roads has not been conducted. 

%FRNDP modeling and solution development is just one of many pre-disaster activities to facilitate better post-disaster recovery efforts. Certain road segments are pre-determined for use by the first responders (after a disaster), and these segments – perhaps just one lane in a multi-lane road - are labeled and communicated clearly to both the citizens and the first responders (FR) in advance. This allows the evacuees a disaster who are mobile and can travel on their own to use only those road segments that are not marked for use by FR, while simultaneously, FR can use these road segments to reach people who are not mobile and require help. 

The specific requirements of where people will be located at the time of the disaster, how many people will be mobile, and likewise, the number of immobile people to whom FRs should reach out, depends on the disaster. However, pre-disaster scenario analysis can provide insights into the impact of various choices of pre-selected FR paths. Every one of such scenarios is a complex mixed-integer non-linear program (MINLP) that is cumbersome to solve, in part because the objective function is itself an outcome of an optimization problem – that of selfish routing of evacuees using paths that are not allocated to a feasible set of FR paths – and not an algebraic function, and we aim to find the best set of FR paths incorporating evacuee behavior constrained by the FR path. 

The primary methodological contribution of this paper is to explore a quantum-inspired approach to solving any one such deterministic instance of FRNDP. We anticipate this will be embedded in an outer loop that studies various scenarios, the parameters of each scenario having been determined by experts in geology and practitioners of humanitarian operations. We hope our work adds to the optimization literature on (a) the disaster preparedness stage and is helpful towards such DM initiatives (not only in Türkiye that specifically motivated this research, but in other regions worldwide as well where such initiatives are being considered) and (b) the emerging area of quantum and quantum-inspired computing that appears promising to tackle complex MINLP models that capture the FRNDP and other complicated, yet critical real-world problems.

\section{Related works}
The \emph{first responder network design problem} (FRNDP) is a variant of the discrete network design problem (DNDP), a well-studied bilevel optimization problem in transportation. The classic DNDP aims to identify an optimal set of candidate links to be adjusted in the network subject to a budget, taking into consideration users' reactions to this alteration, typically governed by a traffic equilibrium perspective \cite{dafermos1980traffic}. The DNDP can be expressed as a bilevel optimization problem, where the leader (outer) problem seeks to optimally determine a network design that minimizes a travel-oriented metric (e.g., total travel time) in the network, and the follower (inner) problem models the users' reactions within the network \cite{bard2013practical}. This reaction is typically represented as a static traffic assignment problem (TAP), taking into account user equilibrium. The first exact method proposed to solve DNDP was a branch and bound method by  \cite{leblanc1975algorithm}, which used the system-optimum relaxation
of the TAP, a convex programming problem, to generate lower bounds. A survey of exact methods to solve the DNDP is given by \cite{rey2020computational}, where it is noted that the largest DNDP instances solved to optimality remain of small scale.

Three primary classes of methods have been proposed to solve the DNDP: a branch-and-bound algorithm \cite{leblanc1975algorithm,farvaresh2013branch,bagloee2017hybrid,yin2022bo}, a generalized Benders' decomposition approach \cite{gao2005solution,fontaine2014benders1,bagloee2017hybrid,fontaine2018benders2} to solve a mixed-integer nonlinear programming formulation, and a single-level formulation created by using the (Karush–Kuhn–Tucker) KKT conditions since the TAP can be represented as a convex nonlinear program \cite{farvaresh2011single,bard2013practical,fontaine2014benders1}. Several other papers solve the DNDP using techniques such as a linear approximation of the travel time function, system-optimum (SO) relaxations of their TAP models (a single-level optimization problem which ignores the follower objective function), and valid inequalities \cite{luathep2011global,wang2013global,wang2015novel}. A hybrid machine learning and bi-objective optimization algorithm is proposed to solve a variant where lane closure and reversal decisions are made at the outer level to mitigate traffic congestion \cite{zhang2023lane}. 

Several variants of DNDP have been proposed to establish emergency routes while considering the evacuation problem at the inner level \cite{afkham2022, nikoo2018, shariat2020}. \cite{afkham2022} employed Benders Decomposition and heuristic accelerators as the solution approach. \cite{nikoo2018} identify the emergency transportation routes that will be closed to public traffic, as in our problem. Three objectives are optimized by a branch and cut algorithm.  \cite{shariat2020} also addresses a multi-objective disaster response routes design problem. However, multiple earthquake scenarios with the probability of failure of each arc are considered in a stochastic programming framework. A branch and cut algorithm is proposed to solve the nonlinear mixed integer program obtained by taking the weighted sum of the objectives. The need to reserve lanes for time-critical operations also arises during sports events such as the Olympics, where the athletes and equipment must travel between the venues within strict time windows. The problem of identifying which lanes to reserve to minimize the total weighted travel time increase for normal traffic while ensuring the time windows was defined as the lane reservation problem by \cite{wu2009}. Ignoring the inner traffic assignment problem in the DNDP, the authors provided a single-level mixed integer linear programming model and a greedy heuristic algorithm for its fast solution. 

The FRNDP introduces additional layers of complexity compared to the classic DNDP. In the DNDP, any feasible combination of arcs in the network can form a solution to the leader problem. However, the set of selected arcs for FR lane reservation in FRNDP must configure paths that connect the nodes requiring FR visits to one of the FR entry points in the network. Due to these constraints, it is not immediately clear how to extend the three methods mentioned above that solve the DNDP. Additionally, the conventional budget constraint significantly restricts the number of feasible solutions in DNDP. In contrast, FRNDP lacks this constraint, differentiating our problem from state-of-the-art variations of DNDP and expanding our search space into a combinatorially large realm. %The complexities mentioned above, as well as the additional FRNDP constraints with respect to the classic DNDP, may render the Benders decomposition \cite{gao2005solution} and KKT methods \cite{farvaresh2011single} less efficient algorithms for FRNDP.

Owing to its non-linear complex objective function subject to linear constraints, our study pioneers in investigating the application of a novel GAGA (GAMA + GAMA) algorithm in solving FRNDP instances. GAGA consists of a bi-level nested GAMA (Graver Augmented Multi-seed Algorithm). GAMA is a Quantum-Classical hybrid algorithm that consists of two main parts: finding the Graver basis and augmenting feasible solutions. GAMA finds successful applications in portfolio management \cite{venturelli2019reverse,alghassi2019graver}, cancer genomics \cite{alghassi2019quantum}, and Quadratic (Semi-) Assignment Problem \cite{alghassi2019GAMA}.

\section{First responder network design problem}
We formally introduce the FRNDP as follows. The FRNDP is defined on a network with nodes $N$ and directed links $A$. Let $F$ be the set of nodes requiring an FR to visit them, and $S$ be the set of nodes with disaster-affected individuals who aim to evacuate the network, referred to as \emph{agents}. Moreover, let $E$ be a set of critical nodes from where FRs can enter the network and agents can leave the network. Representing the travel demand right after the disaster, for each $i \in N$, let $d_i$ be the number of agents (vehicles) at node $i$ who evacuate the network. For each $(i,j) \in A$, let $c_{ij}$ be the capacity of the link, $T_{ij}$ be the free-flow travel time of the link, $l_{ij}$ be the number of lanes on the link.

For each $(i,j) \in A$, let $f_{ij}$ denote the number of agents traveling on link $(i,j)$. For each $(i,j) \in A$, let $t_{ij}(f_{ij})$ denote the travel time from node $i$ to node $j$ using link $(i,j)$, which depends on the flow along the link to account for traffic congestion. We assume that the travel time function is a strictly convex function and has the form $t_{ij}(f_{ij}) = T_{ij}(1 + \alpha (\frac{f_{ij}}{c_{ij}})^{\beta})$, where we set $\alpha = 0.15$ and $\beta = 4$ following the Bureau of Public Roads (BPR) standard in \cite{us1964traffic}. We assume that agents decide to travel to an exit node in $E$, which can be altered depending on the behavior of other agents under optimal selfish routing behavior.

We must decide which lanes to \emph{reserve} for FRs to use. We assume that reserving a lane means reserving one lane in both directions. This assumption is justified by the fact that many FRs will not only need to visit a node but also return to where they came from (and transport immobile people to exit points and so on). So, reserving a lane from a link $(i,j)$ updates the capacities as follows: $c_{ij} \gets \frac{l_{ij} - 1}{l_{ij}}c_{ij}$ and $c_{ji} \gets \frac{l_{ji} - 1}{l_{ji}}c_{ji}$. Note that if more than  one lane exists on a link, the unreserved lanes on this link can be used by mobile evacuees. We also assume that multiple nodes in $F$ can use the same reserved lane to connect to their respective nodes in $E$ without further affecting the capacity. The FRNDP seeks to find a set of lanes to reserve for FRs such that each node in $F$ is connected to a node in $E$ by these reserved lanes, and the sum of the exit times of overall evacuees located nodes in $S$ at user equilibrium traffic is minimized.

\subsection{Mixed integer non-linear programming (MINLP) model}
We model the FRNDP as a mixed integer non-linear program as follows. We will define a bilevel program with outer problem $P_o$ and inner problem $P_i$. The outer problem corresponds to choosing which lanes to reserve for FRs, and the inner problem corresponds to solving the user equilibrium traffic problem (with $\delta_{ik}=1$ if $i=k$ and zero otherwise).

\begin{figure}[ht]
\begin{align}
(P_o) ~~\min\limits_{y} & \displaystyle \sum_{(i,j) \in A} x_{ij} t_{ij}(x_{ij},y_{ij},y_{ji}) \label{F:obj} \\[3ex] 
\textrm{s.t.} & \displaystyle  \sum_{j \in N : (i,j) \in A}{y_{ijk}} - \sum_{j \in N : (j,i) \in A}{y_{jik}} = \delta_{ik} ~& \forall i \in N \backslash E, \ \forall k \in F \label{F:flow-conservation}\\
  & y_{ijk} \leq y_{ij} & \forall (i,j) \in A, \ \forall k \in F \label{F:linking}\\
  & y_{ij} \in \{0,1\} & \forall (i,j) \in A \label{F:binary1}\\
  & y_{ijk} \in \{0,1\} & \forall (i,j) \in A, \ \forall k \in F \label{F:binary2}\\
  & x \in \mathcal{S}_{y}^{*} \label{F:bilevel}
\end{align}

\begin{align}
(P_i) ~~\min\limits_{x} & \displaystyle \sum_{(i,j) \in A}{\int_{0}^{x_{ij}}{t_{ij}(v,y_{ij},y_{ji})dv}} \label{S:obj}\\[3ex]
\textrm{s.t.} & \displaystyle  \sum_{j \in N : (i,j) \in A}{x_{ij}} - \sum_{j \in N : (j,i) \in A}{x_{ji}} = d_{i} ~& \forall i \in N \backslash E \label{S:flow-conservation}\\
  & x_{ij} \geq 0 & \forall (i,j) \in A \label{S:continuous-var}
\end{align}
\end{figure}

The outer problem $P_o$ is defined as follows. For each link $(i,j) \in A$, we define a binary variable $y_{ij}$ that represents the choice of reserving a lane for FRs. For each link $(i,j) \in A$ and each node $k \in F$, we define a binary variable $y_{ijk}$ that represents the choice of an FR using the link to reach node $k$. This second set of binary variables will be useful because of the assumption that multiple nodes in $F$ can be connected to a node in $E$ using the same reserved link. The objective function (\ref{F:obj}) is to minimize the sum of the total travel times of all agents from their initial nodes to an exit. We rewrite the time travel function for a given link $(i,j) \in A$ as $t_{ij}(f_{ij}, y_{ij}, y_{ji}) = T_{ij}\big(1 + \alpha (\frac{f_{ij}}{c_{ij}(y_{ij},y_{ji})})^{\beta}\big)$ to explicitly express that the functions $c_{ij}$ and $t_{ij}$ depend on variables $y_{ij}$ and $y_{ji}$. Constraints (\ref{F:flow-conservation}) ensure that each node in $F$ is connected by reserved lanes to an FR entry node. Constraints (\ref{F:linking}) ensure that a lane from the link $(i,j)$ is considered reserved if one or more nodes in $F$ use the link to connect to an FR entry node. Constraints (\ref{F:binary1} and \ref{F:binary2}) require the $y$ variables to be binary. Constraint (\ref{F:bilevel}) ensures that the solution to $P_o$ uses an optimal solution to $P_i$, where $x$ is the set of inner problem variables and $P_{i}^{*}$ is the set of optimal solutions to $P_i$ for a solution $y$.

Next, we describe the inner problem $P_i$. The model is based on a given feasible solution $y$ to the outer problem $P_o$. For each link $(i,j) \in A$, we define a continuous variable $x_{ij}$ that represents the flow of agents evacuating the network using the link. The objective (\ref{S:obj}) is the standard objective function used in selfish routing problems (\cite{beckmann1956studies}), albeit with the updated travel time function that indicates how link capacities may be updated due to reserved lanes. Constraints (\ref{S:flow-conservation}) ensure that all agents at nodes in $S$ reach an exit node $E$. Constraints (\ref{S:continuous-var}) enforce the non-negativity of the $x$ variables. We assume that allowing a fractional number of agents gives a close enough approximation to enforcing the more realistic requirement that the flow along each link must be an integer number of agents.

%\section{GAMA formulation: FR path and evacuee flow optimization}

\section{Quantum-inspired algorithm}

In this section, we develop a bi-level optimization algorithm for the FRNDP, building on the Graver Augmented Multi-seed Algorithm (GAMA) developed in \cite{alghassi2019graver} and \cite{alghassi2019GAMA}. First, we outline the GAMA algorithm below for the sake of completeness.

\subsection{GAMA}

GAMA is a heuristic algorithm for solving general non-linear integer optimization problems of the form
\begin{equation} \label{eq:IP}
    (\mathcal{IP})_{A, b, l, u, f} = 
            \begin{cases}
            min \ f(x) \\
            Ax = b & l \leq x \leq n \quad x, l, u \in \mathbb{Z}^{n} \\
            A \in \mathbb{Z}^{m \times n} & b \in \mathbb{Z}^{n} \\
            \end{cases},
\end{equation}
where $f: \mathbb{R}^{n} \rightarrow \mathbb{R}$ is a real valued function. One approach to solving such problems is to use an augmentation procedure. An augmentation procedure starts from an initial feasible solution and, at each iteration, checks for an improvement step by using directions from a set of candidate improvement directions, referred to as a test set or optimality certificate. When none of the candidate directions can yield an improvement step, the algorithm terminates. We define the test set below (definition taken verbatim from \cite{alghassi2019graver}).

\begin{definition}
A set $S \in \mathbb{Z}^{n}$ is a test set or an optimality certificate if for every non-optimal, feasible solution, $x_{0}$, there exists $t \in S$ and $\lambda \in Z_{+}$ such that $f(x_{0} + \lambda t) < f(x_{0})$. The vector $t$ is called the \textit{augmenting} direction. 
 \end{definition}

The GAMA algorithm uses a Graver basis \cite{Graver1975OnTF} as a test set. We provide a definition of the Graver basis below (taken verbatim from \cite{alghassi2019graver}). To understand the definition of Graver basis, we first give two auxiliary definitions.

\begin{definition}
Given $x, y \in \mathbb{R}^{n}$, we define x is conformal to y, written as $x \sqsubseteq y$, if $x_{i} y_{i} \geq 0$ (x and y lie in the same orthant),  and $|x_{i}| \leq |y_{i}|$, $\forall$ $i \in \{ 1..n \}$. A sum $u = \sum _{i} v_{i}$ is called conformal if $v_{i} \sqsubseteq u$, $\forall i$.
\end{definition}

For a matrix $A \in \mathbb{Z}^{m \times n}$, define the lattice
\begin{equation} \label{eq:L*}
    L^{*}(A) = \{ x | Ax=0, x \in \mathbb{Z}^{n}, A \in \mathbb{Z}^{m \times n} \} \backslash \{ 0 \}.
\end{equation}

\begin{definition}
    The Graver basis, $\mathcal{G} (A) \subset \mathbb{Z}^{n}$, of an integer matrix $A$ is defined as the finite set of $\sqsubseteq$ minimal elements in $L^{*}(A)$  (i.e, $\forall g_{i}, g_{j} \in \mathcal{G} (A)$, $g_{i}\not \sqsubseteq g_{j}$, when $i \not = j$).
\end{definition}

The Graver basis \cite{Graver1975OnTF} of an integer matrix $A \in \mathbb{Z}^{m \times n}$ is known to be a test set for integer linear programs. Graver basis is also a test set for certain non-linear objective functions, including Separable convex minimization \cite{10.1016/j.orl.2003.11.007}, Convex integer maximization \cite{DELOERA20091569}, Norm $p$ minimization \cite{hemmecke2011polynomial}, Quadratic \cite{{Lee2010TheQG},{10.1016/j.orl.2003.11.007}} and Polynomial minimization \cite{Lee2010TheQG}. It has also been shown that for these problem classes, the number of augmentation steps needed is polynomial \cite{{hemmecke2011polynomial},{DELOERA20091569}}. Graver basis can be computed (only for small size problems) using classical methods such as the algorithms developed by \cite{10.1145/236869.236894} and \cite{Sturmfels1997VariationOC}.

For general non-linear, non-convex optimization problems such as the $\mathcal{IP}$ defined in \ref{eq:IP}, it is not known if the Graver basis is a test set. However, if a non-convex objective function can be viewed as many convex functions stitched together, then the entire feasible solution can be viewed as a collection of parallel subspaces with convex objective functions. Given the Graver basis for the constraint matrix and feasible solutions in each of these subspaces, we can find a global optimal by augmenting along the Graver basis from each of the feasible solutions to reach a local optimal solution and find the best possible solutions among the local optima.

In case we are unsure whether the Graver basis serves as a test set or encounters computational challenges, a pragmatic approach involves leveraging GAMA as a heuristic, utilizing a \emph{partial Graver basis}. A partial Graver basis is a subset of the Graver basis. Using a partial Graver basis, GAMA gives a feasible but not provably optimal solution. Such a method has been studied in \cite{alghassi2019GAMA, alghassi2019graver}, where the heuristic is improved by starting GAMA from multiple initial feasible solutions to the constraint equations as starting points (seeds) instead of just one.
 
\subsection{Partial Graver basis computation}
In general, a partial Graver basis can be obtained as follows. The method finds points in $L^{*}(A)$ by taking the differences of feasible solutions. Then, a classical post-processing step by $\sqsubseteq$-minimal filtering \cite{alghassi2019GAMA} yields a partial Graver basis.

To obtain a sample of feasible solutions to $\mathcal{IP}$, so that we can take their differences, we solve a quadratic unconstrained integer optimization (QUIO) given by
\begin{equation} \label{eq:quio}
\begin{split}
    min \quad x^{T}Q_{I}x - 2b^{T}&Ax \\
    Q_{I} = A^{T}A,
    x \in \mathbb{Z}^{n}.
\end{split}
\end{equation}
If the variable $x$ is an integer, a binary encoding of the integer variables is required. Expressing 
\begin{equation}
    x = L + EX,
\end{equation}
with $L$ as the lower bound vector and $E$ as the encoding matrix, we get a Quadratic Unconstrained Boolean Optimization (QUBO)

\begin{equation}
\begin{split}
    min \quad X^{T}Q_{B}X  \\
X \in \{0, 1 \}^{nk},   
\end{split}
\end{equation}
where $Q_{B} = E^{T}Q_{I}E + 2 diag[(L^{T}Q_{I} -b^{T} A) E]$ and $Q_{I} = A^{T}A$.

The above QUBO can be solved via an annealer or by means of simulated annealing \cite{1983Sci...220..671K} to obtain a sample of feasible solutions. An annealer refers to a specially constructed quantum or semi-classical hardware \cite{{Glover2018QuantumBA},{10.3389/fphy.2014.00005},{mohseni2022ising},{king2018observation},{doi:10.1126/science.aat2025},  {10.1007/978-3-030-19311-9_19},{Chou2019AnalogCO}, {doi:10.1126/science.aah5178},{PhysRevA.88.063853}} used for solving QUBO problems mapped to an Ising model.

%The GAMA heuristic \cite{alghassi2019GAMA, alghassi2019graver} aims to solve the optimization problem, $\mathcal{IP}$ in Equation~\ref{eq:IP} using multiple feasible solutions to the constraint equations as starting points for augmentation and using the partial Graver basis elements as test sets. 

\subsection{GAMA for FRNDP}
We now explain how GAMA attempts to solve FRNDP by only searching on FR paths. The inner problem is solved using a gradient descent algorithm from \cite{leblanc1975efficient}.

 For each node $k \in F$, we obtain a sample of $n_{paths}$ feasible paths, each of which is connected to some FR entry node in $E$. These feasible paths are  obtained by solving a QUBO (implied by equations~\ref{F:flow-conservation}) using a quantum annealing device (such as D-Wave), classical simulated annealing, or by using a path-finding algorithm such as Yen's algorithm for k-shortest paths \cite{yen1971finding}.
 We create feasible solutions for the entire network by combining the feasible FR paths for each point in $F$. (This gives us $M$ seeds.)

%We construct a partial Graver basis in terms of the variables $y_{ijk}$ to search over a subset of feasible solutions for these variables. Then, at each iteration, the heuristic will solve 

%We create a partial Graver basis to search over feasible solutions of $y_{ijk}$ as follows. We consider the partial Graver basis of a matrix $A_{FR}$, which is defined as a matrix containing only constraints \ref{F:flow-conservation} from $P_o$.

A set of lattice elements (recall equation \ref{eq:L*}) can be obtained  by taking differences of these feasible solutions. Then, a partial Graver basis is obtained from these lattice elements using the procedure in \cite{alghassi2019GAMA, alghassi2019graver} that ensures that all elements are conformal minimal.
%for each node $k \in F$ using the difference of feasible paths as elements of $L^{*}(A_{FR})$ and classical post-processing to obtain conformally minimal (indecomposable) Graver elements. We refer the reader to \cite{alghassi2019GAMA, alghassi2019graver} for detailed numerical examples on generating Graver basis. 

%We combine all of these elements to obtain a partial Graver basis to be used by GAMA.

As in \cite{alghassi2019GAMA, alghassi2019graver}, we start the augmentation algorithm from $M$ initial feasible solutions. Rather than take all feasible solutions as seeds, we construct a random sample of $M$ feasible solutions by combining randomly chosen paths for each $k \in F$ from the corresponding set of $n_{paths}$ paths that we obtained when calculating the partial Graver basis. So, for each initial feasible solution, we randomly select a path for each $k \in F$, and combine the paths to obtain a feasible solution of $y_{ijk}$ variables. An assignment of the $y_{ij}$ variables follows naturally. %A feasible solution to the FR path can be obtained as a linear summation of one FR path to each node $k \in \mathrm{F}$. %Here, the FR path to each node is chosen randomly from the $n_{paths}$ sample of feasible paths. %To do this, for each link $(i,j) \in A$, let $z_{i,j}^{r}$ be an integer variable for a random solution instance, $r \in \{1,...,M\}$. \ananth{elaborated explanation now} Then,

%\begin{equation} \label{F_random_sol}
%    z_{ij}^{r} = \sum _{k \in F} y_{ijk}^{r_{k}}
%\end{equation}
%where, $y_{ijk}^{r_{k}}$ is a binary variable and $r_{k} \in \{1,..,n_{paths}\}$ is a random integer.

%Now, the binary variable, $y_{ij}=1$ (indicating if a link $(i,j) \in A$ is reserved for FRs) if $z_{ij} \geq 1$. The partial Graver basis for the FR paths can now be obtained as a concatenation of the partial Graver basis elements for FR paths to each node in $F$. For a given FR feasible solution, the objective function is the total evacuation time (\ref{F:obj}), subject to the road network capacity under the given FR and the evacuee flow, $x_{ij}$, which is obtained as an equilibrium flow solution under user equilibrium assumption. The evacuee flow is obtained by optimizing the user equilibrium objective function (\ref{S:obj}) in the inner problem, $P_i$,  with a fixed FR path. The inner problem can be solved using a gradient descent algorithm as in ~\cite{leblanc1975efficient}. However, we noticed that the inner problem itself can also be solved using GAMA, which improves the running time of the algorithm. We explain the GAMA approach to solve the inner problem $P_i$ below.

Using a partial Graver basis and $M$ initial seeds, GAMA works as follows. For each initial seed, start from the corresponding initial feasible solution. Then, check each direction in the partial Graver basis to update the variables $y_{ijk}$ and $y_{ij}$. To check a direction, use step size $t=1$ and run the gradient descent algorithm from \cite{leblanc1975efficient} to calculate the values of the $x$ variables. If a step is improving, the solution is updated along the step and the directions in the Graver basis are checked again. The algorithm terminates when no step gives an improved objective function value from the current one. 

We illustrate GAMA using an example.

%We provide an illustrative example using a $4-node$ graph. There are four FR paths in all, and 14 Graver elements. 
\subsubsection{Example with a 4-Node Network}
Consider a $4$-node network graph, $G = (N,A)$ with $N = \{ 0, 1, 2, 3 \}$ and $A = \{ (0,1), (0,2), (0,3), (1,2),\\ (1,3), (2,3) \}$ with one lane in each edge. If a FR lane is allocated to an edge, no demand can flow through it in the opposite direction. For simplicity, we will also assume a single entry/exit: node $3$. The demand is also present at a single location: node $0$. So, we have set of nodes requiring FR, $F = \{ 0\}$ and critical nodes, $E = \{ 3\}$ and agent node, $S = \{ 0\}$. Let $d_{0}=100$. The arrows in the graph show the evacuee direction via different edges, if available, from node $0$ to node $3$.  Recall that there is only one lane for each edge. Therefore, assigning an FR path to an edge makes it unavailable for evacuating. The capacities on the edges are  $C = \{25, 30, 35, 35, 15, 45 \}$, and the free flow travel time is taken as 1. \\

\begin{figure}
\begin{center}
    \begin{tikzpicture}   
  \node [shape=circle,draw=black] (a0) at (2, 3) {$0 _{\ d_{0}=100}$};  
  \node [shape=circle,draw=black] (a1) at (4, 5)  {$1_{\ d_{1} = 0}$};
  \node [shape=circle,draw=black] (a2) at (4, 1)  {$2_{\ d_{2} = 0}$};  
  \node [shape=circle,draw=black] (a3) at (8, 3) {$3_ {entry/exit}$}; 
  
  \draw [arrow] (a0) -- (a1) node[midway, right]{};
  \draw [arrow] (a0) -- (a2) node[midway, right]{};  
  \draw [arrow] (a0) -- (a3) node[midway, above]{};  
  \draw [arrow] (a1) -- (a2) node[midway, below]{};  
  \draw [arrow] (a1) -- (a3) node[midway, right]{};  
  \draw [arrow] (a2) -- (a3) node[midway, right]{};   
    \end{tikzpicture}  
\end{center}
\caption{An example of a $4$-node network.}
\end{figure}
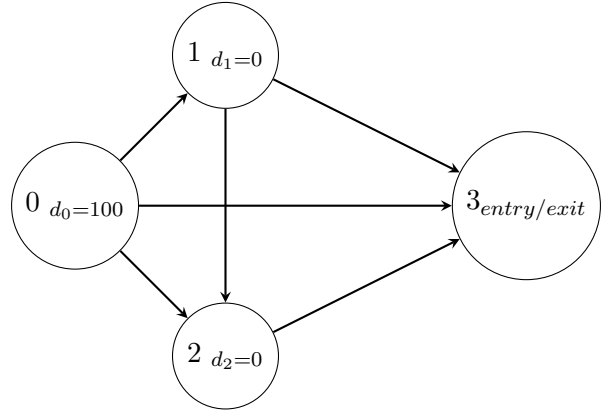

We suppress the index $k=0$ here as it is the only demand node. With $y = (y_{01}, y_{02}, y_{03},  y_{12}, y_{13}, y_{23})$, the FR path constraints are

\begin{equation}
     y_{01} + y_{02} + y_{03} = 1 
\end{equation}
\begin{equation}
    y_{12} + y_{13} - y_{01}= 0 
\end{equation}
\begin{equation}
    y_{23} - y_{02} - y_{12} = 0 
\end{equation}
where $y_{ij} \in \{ 0, 1 \}$.

The above constraints can be expressed as $A_{FR} X_{FR} = B_{FR}$ with 
\begin{center}
$X_{FR} = \begin{pmatrix} 
    y_{01} \\
    y_{02} \\ 
    y_{03} \\ 
    y_{12} \\ 
    y_{13} \\ 
    y_{23}
\end{pmatrix}$,   
$A_{FR} = \begin{pmatrix} 
     1 & 1  & 1 &  0 & 0 & 0 \\
    -1 & 0  & 0 &  1 & 1 & 0 \\ 
     0 & -1 & 0 & -1 & 0 & 0 
\end{pmatrix}$,
$B_{FR} = \begin{pmatrix} 
     1\\
     0\\
     0
\end{pmatrix}.$
\end{center}

There are a total of 4 feasible solutions given by
 \begin{center}
$Y_{FR, feasible} = \begin{pmatrix} 
    1 \\ 0 \\ 0 \\ 1 \\ 0 \\ 1
    \end{pmatrix},
    \begin{pmatrix} 
    1 \\ 0 \\ 0 \\ 0 \\ 1 \\ 0
    \end{pmatrix},
    \begin{pmatrix} 
    0 \\ 1 \\ 0 \\ 0 \\ 0 \\ 1
    \end{pmatrix},
    \begin{pmatrix} 
    0 \\ 0 \\ 1 \\ 0 \\ 0 \\ 1
    \end{pmatrix}.
    $ 
\end{center}

The lattice elements satisfying $A_{FR} X_{FR} = 0$ can be obtained as a difference of the feasible solutions given by
\begin{center}
$Y_{FR, kernel} = \begin{pmatrix} 
    0 \\ 0 \\ 0 \\ 1 \\ -1 \\ 1
    \end{pmatrix},
    \begin{pmatrix} 
    1 \\ -1 \\ 0 \\ 1 \\ 0 \\ 0
    \end{pmatrix},
    \begin{pmatrix} 
    1 \\  0 \\ -1\\  1\\  0\\  1
    \end{pmatrix},
    \begin{pmatrix} 
    1\\ -1\\  0\\  0\\  1\\ -1
    \end{pmatrix}
    \begin{pmatrix} 
    1\\  0\\ -1\\  0\\  1\\  0
    \end{pmatrix}
    \begin{pmatrix} 
    0\\  1\\ -1\\  0\\  0\\  1
    \end{pmatrix}.
    $
\end{center}
The lattice elements are already conformally minimal elements and so form a subset of the full Graver basis. So, the lattice elements $Y_{FR, kernel}$ is a partial Graver basis set in this instance. The full Graver basis (obtained using Pottier's algorithm \cite{10.1145/236869.236894}) is given by
\begin{center}

$\mathcal{G}(A_{FR})=
\begin{pmatrix}
    1  &  1 &  0 & -1 & -1 &  0 &  1 & -1 &  0 &  0 &  1 & -1 &  0 &  0 \\
    -1 &  0 &  1 &  1 &  0 & -1 &  0 &  0 & -1 &  1 & -1 &  1 &  0 &  0 \\
    0  & -1 & -1 &  0 &  1 &  1 & -1 &  1 &  1 & -1 &  0 &  0 &  0 &  0 \\
    1  &  0 &  0 & -1 &  0 &  0 &  1 & -1 &  1 & -1 &  0 &  0 &  1 & -1 \\
    0  &  1 &  0 &  0 & -1 &  0 &  0 &  0 & -1 &  1 &  1 & -1 & -1 &  1 \\
    0  &  0 &  1 &  0 &  0 & -1 &  1 & -1 &  0 &  0 & -1 &  1 &  1 & -1  \\ 
\end{pmatrix}$.
\end{center}

Given how small the example is, we can enumerate all FR path choices and compute the resulting total evacuation time in equilibrium given that FR choice. See Table~\ref{tab:4-node_fr}. Note the almost identical travel time for each evacuation path under user equilibrium at each FR choice. We can see that FR Path - 2 is optimal for the objective of total time for evacuation.  
\begin{table}
\centering
    \begin{tabular}{|c|c|c|c|c|c}
    \hline
        ID & FR-Path & Evac. paths & Path evac. time & Total evac. time\\
    \hline
        $1$ & $(3,0)$ & $P_{1}:(0, 1, 3)$    & 5.02  & 507.56 \\  
            &         & $P_{2}:(0, 1, 2, 3)$ & 5.18  &        \\ 
            &         & $P_{3}:(0, 2, 3)$    & 5.08  &        \\ 
    \hline
        $2$ & $(3,1,0)$ & $P_{1}:(0, 2, 3)$  & 2.51  & 243.43\\
            &           & $P_{2}:(0,3)$      & 2.38  &       \\
    \hline
        $3$ & $(3,2,0)$ & $P_{1}: (0, 1, 3)$ & 4.06  & 379.01\\
            &           & $P_{2}: (0, 3)$    & 3.67  &       \\
    \hline
        $4$ & $(3,2,1,0)$ & $P_{1}: (0, 3)$  & 11.00 & 1099.58\\
    \hline
    \end{tabular}
    \caption{Agent evacuation times for each FR in a $4-$node example.}
    \label{tab:4-node_fr}
\end{table}

Now, we will illustrate steps in GAMA. Suppose we use only one of the feasible FR solution as a seed (M=1), namely $(3,0)$ with ID-1 in Table~\ref{tab:4-node_fr}. The solution is $X_{FR} = (0 \ 0 \ 1 \ 0 \ 0 \ 0 )^{T}$. We solve the inner problem, $\mathcal{P}_{i}$ to obtain the optimal user equilibrium objective. The optimal flow values are found to be $X_{FL}^{*} = (42 \ 58 \ 0 \ 14 \ 28 \ 72 )^{T}$. Using these flow values on the edges, we evaluate the total evacuation time, which is $T = 507.56$. Note that the travel time along the $3$ possible evacuation routes, $P_{1}, P_{2}, P_{3}$ are approximately similar. They are not the same because of integrality constraints.

Now, we start the Graver augmentation. Consider the direction with Graver element, $g_{0} = (1 \ -1 \ 0 \ 1 \ 0 \ 0 )^{T}$. Clearly, $X_{FR}^{0} + g_{0}$ is infeasible. Consider the next Graver element, $g_{1} = (1 \ 0 \ -1 \ 0 \ 1 \ 0 )^{T}$. Here $X_{FR}^{1} = X_{FR}^{0} + g_{1} = (1 \ 0 \ 0 \ 0 \ 1 \ 0 )^{T}$ which is the FR-path with ID-2. Solving the inner problem in this case gives the optimal flow values, $X_{FL}^{*} = (0 \ 39 \ 61 \ 0 \ 0 \ 39)^{T}$ with a total evacuation time of $T=243.43$. This is better than $T=507.56$. So, we update the optimal FR path to $X_{FR}^{1}$.

Continuing with the Graver walk, the next feasible solution is $X_{FR}^{1} + g_{4}$, with $g_{4} = (-1 \ 0 \ 1 \ 0 \ -1 \ 0)^{T}$ which is again $X_{FR}^{0}$. Since the evacuation time is not improved, the FR path remains current at $X_{FR}^{1}$. The next feasible FR path in the course of the Graver walk is $X_{FR}^{3} = X_{FR}^{2} + g_{11}$ with $g_{11} = (-1 \ 1 \ 0 \ 0 \ -1 \ 1)^{T}$ and $X_{FR}^{3} = (0 \ 1 \ 0 \ 0 \ 0 \ 1)^{T}$, $X_{FL}^{*} = (35 \ 0 \ 65 \ 0 \ 35 \ 0)^{T}$ with evacuation time of 379.01 and then $X_{FR}^{4} = (1 \ 0 \ 0 \ 1 \ 0 \ 1)^{T}$, $X_{FL}^{*} = (0 \ 0 \ 100 \ 0 \ 0 \ 0 )^{T}$ using the Graver element, $g_{12} = (0 \ 0 \ 0 \ 1 \ -1 \ 1 )^{T}$ with evacuation time of 1099.58. At this point, we have verified that no further augmentation is possible from $X_{FR}^{2}$ and the algorithm terminates with $X_{FR}^{2}$ as the optimal FR path solution.

\subsection{GAGA: Bi-level nested GAMA within GAMA} \label{S:gaga}
A further methodological enhancement of this work is to attempt to solve the inner optimization with GAMA itself. The motivation is that approximately solving $P_i$ at each step will accelerate the search, because exactly solving $P_i$ has been difficult.  The solution to $P_i$ that GAMA finds, however, is an approximation because GAMA uses a partial Graver basis for its search. However, upon termination, our procedure will add another step: solve $P_i$ using the exact solution method by \cite{leblanc1975efficient}.
To the best of our knowledge, this is the first time that a heuristic like this has been studied. We call this algorithm that is GAMA within GAMA as GAGA.

%(discussed in Section~\ref{S:BB} below) to achieve an optimal solution, ensuring that the heuristic returns a feasible solution to $P_o$. 
%We develop a bi-level Graver augmented optimization approach to solve the FRNDP by solving both the outer problem, $P_o$ and the inner user equilibrium problem, $P_i$ using the GAMA heuristic. The method to solve $P_o$ using GAMA has been discussed in the previous subsection.

We construct each initial feasible solution as follows. For each node $k \in S$, we randomly and independently sample $d_{k}$ paths from the set of $n_{paths}$ paths and send one evacuee along each path. As before, to generate a partial Graver basis with respect to the constraint matrix in $P_i$, we find lattice elements that are differences of paths between nodes with evacuees in $S$ and exit nodes in $E$, and then apply the same procedure to end up with conformally minimal solutions. To create many of these paths, for each node, $k \in S$, we generate a sample of $n_{paths}$ feasible paths to exit nodes in $E$. As before, this can be done with an annealer, simulated annealing, or a path-finding algorithm.

The inner GAMA solves $P_i$ as follows. We start from each of the $M$ different initial solutions. At each step, the algorithm checks each possible improving direction in the Graver basis, using a step size of $t=1$, and updates the current solution as soon as an objective improving direction is found. We terminate the inner GAMA when there is no improving direction up to a tolerance. Without a tolerance threshold, the Graver augmentation runs until there is no improving direction amongst all the Graver basis elements. A tolerance threshold ensures that the run is assumed to be completed if the change in the objective function after one complete pass over all the Graver basis elements is below the tolerance threshold. We experiment with and without a tolerance threshold of $10^{-3}$ for the inner problem.

We note that the Graver augmentation for $P_i$ is carried out using both a partial Graver basis and integral flow values. Hence, it is possible and can be expected that the user equilibrium flow obtained by the GAMA heuristic is a close approximation to, but not the exact user equilibrium flow solution. Therefore, as a final step, we run the gradient descent algorithm from ~\cite{leblanc1975efficient} on each of the local optimal FR solutions at the end of the GAGA and report the optimal total evacuation time under this user flow equilibrium as the optimized result. In the rest of the paper, we denote the final solution obtained at the end of the GAGA procedure as GAGA-\textit{only} and the solution after running the additional step of the gradient descent method from~\cite{leblanc1975efficient} as GAGA-\textit{LeBlanc}. Unless explicitly stated otherwise, the final solutions presented refer to the GAGA-\textit{LeBlanc} values.

To summarize, the GAGA algorithm works as follows. We run the (outer) GAMA heuristic to solve $P_o$ and at each step we solve $P_i$ using (inner) GAMA up to a tolerance. In order to speed up the computation, we use the current equilibrium flow solution as the initial feasible solution to the inner problem. The GAGA algorithm terminates when none of the directions in the partial Graver basis (outer and inner) can improve the current solution.

%All the FR path choices and the resulting total evacuation times in equilibrium, given that FR choice are enumerated in Table~\ref{tab:4-node_fr}

\section{Branch-and-Bound algorithm} \label{S:BB}
%We develop a branch-and-bound algorithm to solve the FRNDP in order to better evaluate the performance of GAGA. The method is based on an algorithm that was developed by~\cite{leblanc1975algorithm} to solve the DNDP and later extended by . The algorithm is adapted to consider the additional constraints of the FRNDP regarding reserving links for first responders. Next, we will describe the details of the branch-and-bound algorithm including the primal heuristic, branching rule, and dual-bound heuristic. We recognize that more sophisticated heuristics can be developed and tested, but we believe that this version of branch-and-bound is a worthy comparison for GAGA.

We develop a branch-and-bound algorithm to solve the FRNDP to evaluate the performance of GAGA. The method is based on an algorithm that was first developed by~\cite{leblanc1975algorithm} and later extended in \cite{farvaresh2013branch} by solving the system optimal relaxation of the DNDP—instead of the traffic assignment problem. The algorithm is adapted to consider the additional constraints of the FRNDP regarding reserving links for first responders. Next, we will describe the details of the branch-and-bound algorithm, including the primal heuristic, branching rule, dual-bound heuristic, and node selection heuristic. We recognize that more sophisticated heuristics can be developed and tested. Still, we believe that this version of the branch-and-bound algorithm which is considered among the most efficient exact solution methodologies in the literature \cite{rey2020computational} provides a worthy comparison for GAGA.

We use the following primal heuristic to find a feasible solution to $P_o$. First, we compute the shortest path from each node in $E$ to each node in $F$ using Dijkstra's algorithm. Then, for each node in $F$, we reserve the links along the shortest of these paths, which corresponds to an assignment of the variables $y$. Next, we solve $P_i$ by using a gradient descent algorithm from~\cite{leblanc1975efficient}. This gradient descent algorithm solves each subproblem by several shortest path calculations for each node in $S$. The algorithm uses a halving method to calculate the step size, and we use a relative tolerance of $0.001$ to determine when to terminate gradient descent, which we denote as $\epsilon$. We use a greedy heuristic at the root node to find an initial feasible solution to begin the gradient descent algorithm. At a node that is not the root node, to accelerate the gradient descent algorithm, we initialize a feasible solution with the optimal solution of the parent node.

Next, we describe the branching rule. %This differs from the branching rule in~\cite{leblanc1975algorithm} due to the additional constraints in the FRNDP.
We develop a branching procedure that creates two child nodes in the branch-and-bound tree based on a given link: one that requires a solution to reserve this link and one that requires a solution not to reserve this link. So, we develop a method of choosing the particular link to use in the branching procedure at a given node. We maintain an order of the nodes in $F$, so given the link that created the branch for the current node based on some $i \in F$, we consider the next node $j \in F$ based on the order. Given this node, based on the solution $y$, we order the links that connect $j$ to an exit node $k \in E$, starting from the exit node. We select the first link that has not yet been fixed to be used or not used at the current node.

We can calculate a dual bound by solving for the total evacuation time of a \emph{system optimal} (i.e., authoritarian) solution. We solve this problem by using the same gradient descent algorithm as the one described above. Similar to previous literature, we notice that the lower bounds are not effective in fathoming search nodes in the branch-and-bound tree. Thus, to improve the algorithm's efficiency, we forego calculating lower bounds in our experiments. It would likely be useful to implement lower-bound computations similar to the ones done in \cite{farvaresh2013branch}. As a consequence of not calculating lower bounds, we implement a breadth-first node selection heuristic instead of another commonly used heuristic that chooses the next node to solve based on lower bounds.

\section{Experimental results}
We evaluate the performance of GAGA on two types of instances: random instances and instances created with data from Istanbul, Türkiye, which is under serious earthquake risk. As pointed out in recent news \cite{Istanbul2023}, Istanbul is expected to be hit by an earthquake with a magnitude above 7.0 in the upcoming seven years at a probability of 64 percent. Therefore, the problem we solve carries great significance for the preparedness of Istanbul.  

In our experiments, we compare the performance of GAGA to that of the branch-and-bound method, which we will refer to as BB. Throughout the experiments, we vary the following parameters of GAGA. 
\begin{itemize}
    \item The number of reads or samples in the simulated annealing step or using a D-Wave annealer. We denote this value as $n_{samples}$ and compare $n_{samples} = 10000$ and with $n_{samples}=1000$.
    \item The number of paths from/to each node $i \in N\backslash E$. We set $n_{paths}=25$ as a default and experiment with other values of $n_{paths}$ for a single instance.
\end{itemize}

\subsection{Random graph instances} \label{sec: random graphs}
We create the first set of instances by constructing small random graphs and randomly generating the other relevant data. These graphs serve as an initial evaluation of GAGA's performance and a comparison of its performance with that of BB. 

\subsubsection{Instance generation}
We create randomly generated instances as follows. We start with input parameters for the number of nodes $n$ and the probability of an edge being between two nodes $p$. We produce $n$ nodes on a two-dimensional grid, randomly generating each coordinate from a uniform [0,1] distribution. Then, we produce $\frac{n}{10}$ nodes as exits nodes $E$ on the boundary of the unit square. Then, for each pair of nodes, we create edges in the graph with a probability that is the product of $p$ and an additional factor based on the distance between the two nodes. The capacities of each edge are generated from a normal distribution with an average of $50$ and a standard deviation of $20$. We assign each edge to have two lanes in both directions. We create the demands $d$ at each node from a normal distribution with mean $100$ and standard deviation $10$. We generate ten instances of each type and denote the instance number by $i$. If an instance solution is found to be infeasible, we create another instance until we have ten instances with feasible solutions. The random instances are generated with the number of nodes $n=10, 20, 30$, each with $p=0.5$ and $p= 0.75$.

\subsection{Details of processor used for computational experiments}

The GAGA is implemented on a single core of Intel Xeon Gold 6252 @ 2.10 GHz (Processor-G) with 192GB memory and 48 cores. The Branch and bound have been carried out on a single core of Intel(R) Xeon(R) Gold 6248R CPU @ 3.00GHz (Processor-B) with 64 GB memory and 8 cores. Based on a benchmarking comparison in \textit{PassMark}, the single-core performance of Processor-B is found to be better. Since no numerical factor is provided there, we compared the run times between the two processors for the $n=30, p=0.75$ case. We found (see Appendix A) that Processor-B is faster in all the $10$ random instances and by a factor of $\sim 1.3$ on average. Regardless, when we report the times, we have not scaled it by any factor.\footnote{We further note that only a single core is used for the computational experiments in both the GAGA and Branch-and-Bound runs. There is no parallel processing involved, although there appears good possibility of developing such algorithms. So, the different number of cores in the two processors is not relevant for the comparison of run times.}

\subsubsection{Testing the performance of D-Wave}

Using random graphs, we test the applicability and performance of the D-Wave\footnote{\url{https://www.dwavesys.com/}} quantum annealing device in generating a sample of feasible solutions required for the GAGA algorithm. We generate the feasible solutions for $n=10$ random instances with $p=0.5, 0.75$. The experiments have been performed on the D-Wave Advantage-system4.1 machine with an annealing time of $1 \mu s$ and $1000$ samples. Following the annealing schedule, a post-processing step using a steepest descent solver from D-Wave Ocean tools is carried out. This step is found to be necessary to obtain a larger sample of feasible solutions. The feasible solutions to the constraint equations in Equations~\eqref{F:flow-conservation}, ~\eqref{S:flow-conservation} generate a set of feasible paths from a given node, $i \in N\backslash E$ to one of the critical nodes, $E$. However, the constraints, as shown in Equations~\eqref{F:flow-conservation}, ~\eqref{S:flow-conservation}, also allow cyclic paths as part of the solutions. We run a Dijkstra algorithm starting from the demand node, $k \in N \backslash E$, to one of the exit nodes, using only the arcs with $y_{ijk}=1$ in the feasible solution. For dense graphs, we found that the cyclic paths are obtained at the nodes present along the path from the node $k \in N \backslash E$ to one of the exits in $E$. This can be eliminated by adding the following additional constraint system to the QUBO, requiring that there is only a maximum of one outward edge at each node.  

\begin{equation} \label{S:constraint_1edge}
        \sum _{j \in N: (i,j)\in A} y_{ijk} + a_{ik} = 1, \forall i \in N \backslash E 
\end{equation}
where $a_{ik} \in \{0, 1 \}$ is an auxiliary variable for a given node $k$, from which a feasible path to the critical nodes $E$ is needed.

After generating the feasible solutions and post-processing, we compute the partial Graver basis and carry out a Graver augmentation, followed by a final \cite{leblanc1975efficient} algorithm call, as alluded to in Section~\ref{S:gaga}. We similarly generate feasible solutions using simulated annealing. We use the D-Wave Ocean tools' classical Neal sampler for carrying out the simulated annealing with $1000$ and $10000$ samples. The results are then compared with the Branch-and-Bound (BB) algorithm. We tabulate the results along with the total computational time in Table~\ref{tab:random_instances_N10_Dwave}. The BB algorithm is given a time of $600$ seconds for this instance as the total time with the GAGA algorithm in all the cases is below this time. As seen from the table, the GAGA solutions obtained using D-Wave solutions (DWaveAdv4.1-1000) and the simulated annealing solutions with  $n_{samples} =1000, 10000$ reads (GAGA-SA-1000, GAGA-SA-10000) are comparable for both $p=0.5$ and $p=0.75$. Compared with Branch-and-Bound (BB), it can be seen that GAGA results are similar for $p=0.5$ and are relatively better at $p=0.75$.  

The D-Wave annealer generates feasible solutions only for the graphs of node size $n=10$. We tested the D-Wave annealer using one instance of a graph with $n=20$ nodes ($p=0.5, i=0$) with both the Advantage-System 4.1 and Advantage-System 6.3 machines using $n=1000$ samples. However, we cannot generate feasible solutions for all the demand nodes. Hence, we only present results with feasible solutions generated using simulated annealing in the rest of the paper.

\begin{table}[hbt!]
\fontsize{10}{10}\selectfont
\renewcommand{\arraystretch}{1.2}
\centering
%\begin{adjustwidth}{-1.5cm}{}
    \begin{tabular}{|c|c|c|c|c|c|c|c|c|c|c|}
    \hline
        \multicolumn{3}{|c|}{Instance} & \multicolumn{2}{|c|}{GAGA-SA-1000} & \multicolumn{2}{|c|}{DWaveAdv4.1-1000} & \multicolumn{2}{|c|}{GAGA-SA-10000} & \multicolumn{2}{|c|}{Branch-And-Bound} \\
    \hline
        n & p & i & obj & time (s) & obj & time (s) & obj & time (s) & obj & time (s) \\ \hline
    \hline
        10 & 0.5 & 0 &  6.59~e+07  &  19.81  &  \textbf{6.42~e+07}  &  21.33  &  6.65~e+07  &  86.70  & 6.60~e+07 & 600 \\
        10 & 0.5 & 1 &  \textbf{8.08~e+06}  &  161.53  &  8.49~e+06  &  118.71  &  8.18~e+06  &  251.63  & 8.16~e+06 & 600 \\
        10 & 0.5 & 2 &  1.27~e+08  &  45.98  &  \textbf{1.26~e+08}  &  43.37  &  1.26~e+08  &  123.12  & 1.27~e+08 & 600 \\
        10 & 0.5 & 3 &  11.52~e+08  &  33.43  &  \textbf{4.15~e+08}  &  33.50  &  4.15~e+08  &  101.85  & 4.16~e+08 & 600 \\   
        10 & 0.5 & 4 &  \textbf{5.52~e+07}  &  94.43  &  5.53~e+07  &  73.67  &  5.53~e+07  &  176.50  & 5.53~e+07 & 600 \\
        10 & 0.5 & 5 &  \textbf{2.20~e+08}  &  140.90  &  2.24~e+08  &  123.16  &  2.21~e+08  &  234.15  & 2.21~e+08 & 600 \\
        10 & 0.5 & 6 &  \textbf{6.67~e+06}  &  150.72  &  8.06~e+06  &  162.21  &  6.93~e+06  &  247.18  & 6.80~e+06 & 600 \\
        10 & 0.5 & 7 &  9.12~e+07  &  51.76  &  \textbf{1.23~e+07}  &  52.02  &  2.60~e+07  &  115.27  & 6.81~e+07 & 600 \\
        10 & 0.5 & 8 &  4.24~e+06  &  48.42  &  \textbf{4.08~e+06}  &  47.21  &  4.44~e+06  &  125.42  & 4.13~e+06 & 600 \\
        10 & 0.5 & 9 &  6.79~e+07  &  177.53  &  7.85~e+07  &  138.79  &  \textbf{6.78~e+07}  &  276.38  & 1.01~e+08 & 600 \\
    \hline
        10 & 0.75 & 0 &  1.99~e+05  &  241.60  &  2.08~e+05  &  241.32  &  \textbf{1.97~e+05}  &  413.56  & 2.29~e+05 & 600 \\
        10 & 0.75 & 1 &  5.06~e+05  &  276.44  &  5.54~e+05  &  239.37  &  \textbf{4.96~e+05}  &  405.03  & 6.27~e+05 & 600 \\
        10 & 0.75 & 2 &  \textbf{4.84~e+06}  &  182.41  &  5.36~e+06  &  172.65  &  4.85~e+06  &  285.60  & 6.89~e+06 & 600  \\
        10 & 0.75 & 3 &  3.05~e+05  &  258.83  &  3.10~e+05  &  188.40  &  \textbf{2.91~e+05}  &  453.60  & 3.32~e+05 &  600\\
        10 & 0.75 & 4 &  1.88~e+06  &  290.05  &  \textbf{1.74~e+06}  &  239.80  &  1.86~e+06  &  445.27  & 1.98~e+06 & 600 \\
        10 & 0.75 & 5 &  3.78~e+05  &  233.02  &  \textbf{3.76~e+05}  &  145.28  &  3.82~e+05  &  433.21  & 5.81~e+05 & 600 \\
        10 & 0.75 & 6 &  \textbf{1.81~e+05}  &  227.59  &  1.92~e+05  &  248.97  &  1.81~e+05  &  418.50  & 2.14~e+05 & 600 \\
        10 & 0.75 & 7 &  7.70~e+04  & 263.95 &  \textbf{7.20~e+04}  &  263.13  &  7.70~e+04  &  452.06  & 7.67~e+04 & 600 \\
        10 & 0.75 & 8 & \textbf{2.69~e+06} &  292.89  &  2.69~e+06  &  354.23  &  2.69~e+06  &  446.26  & 5.07~e+06 & 600  \\
        10 & 0.75 & 9 &  \textbf{1.61~e+05}  &  316.29  &  1.71~e+05  &  214.42  &  1.61~e+05  &  496.67  & 1.74~e+05 & 600 \\
    \hline
    \end{tabular}
    %\end{adjustwidth}
    \caption{Results comparing D-Wave Advantage 4.1 output using 1000 samples (with an anneal time of $1\mu s$ per sample) and simulated annealing using 1000, 10000 samples on random instances with n=10 nodes randomly generated in U(0,1), p probability of an edge with an additional factor based on distance, $i$ instance number. An additional $\frac{n}{10}$ nodes are marked as exits on the boundary of the unit square. These also serve as FR entry points. Capacities of edges are N(50,20). All edges have 2 lanes in both directions. The population at each node to exit is N(100,10). For each instance, the entry for the method with the best performance is in bold text.}

    %We ran D-Wave Advantage 4.1 and D-Wave Advantage 6.1 on a single $n=20, p=0.5, i=0$ instance and found that the feasible paths are not found for all the demand nodes in the instance. So, we limit the application of D-Wave annealer only to problem sizes of $n=10$ nodes.
    \label{tab:random_instances_N10_Dwave}
    
\end{table}

\subsubsection{Path generation methods}

In Table~\ref{tab:random_instances_comparison}, we compare the GAGA algorithm implemented using feasible solutions generated with simulated annealing and the Branch-and-Bound results. Additionally, we generate a set of feasible solutions using the Yen's $k$-shortest path algorithm as an alternative to simulated annealing. In all these instances, the BB algorithm is given a run time of $14000 s$, which is longer than the maximum time taken for the GAGA run on any of these instances. The simulated annealing runs are carried out with $n_{samples} = 1000$ and $10000$ to test the dependence of the quality of the solution on the $n_{samples}$ parameter. We refer to these as GAGA-SA-1000 and GAGA-SA-10000, respectively. The results using feasible solutions from Yen's algorithm are referred to as GAGA-Yens. From the table, we can see that the solutions are comparable among GAGA-SA-1000, GAGA-SA-10000, as well as GAGA-Yens. Compared with Branch-and-Bound, we can see that the results from the GAGA algorithm have a better solution quality. We observed here (not shown in any table for space reasons) that using $n_{samples} = 10000$, the solution obtained with the GAGA algorithm before the \cite{leblanc1975efficient} algorithm call (i.e., GAGA-\textit{only} results) is better than the solutions for the same obtained using $n_{samples}=1000$ or the Yen's algorithm. We attribute this to a lack of diversity on the partial Graver basis. We explore this in detail for a single instance in the subsection below.

\begin{table}
\centering
\begin{adjustwidth}{-0cm}{}
\fontsize{10}{10}\selectfont
\renewcommand{\arraystretch}{1.1}
    \begin{tabular}{|c|c|c|c|c|c|c|c|c|c|c|}
    \hline
        \multicolumn{3}{|c|}{Instance} & \multicolumn{2}{|c|}{GAGA-SA-1000} & \multicolumn{2}{|c|}{GAGA-SA-10000} & \multicolumn{2}{|c|}{GAGA-Yens} & \multicolumn{2}{|c|}{Branch-And-Bound} \\
    \hline
        n & p & i & obj & time (s) & obj & time (s) & obj & time (s) & obj & time (s) \\ \hline
    \hline

        20 & 0.5 & 0 &  2.83~e+06  &  2931.32  &  \textbf{2.64~e+06}  &  4368.71  &  2.69~e+06  &  3697.60  & 4.03~e+06 &  14400 \\
        20 & 0.5 & 1 &  1.45~e+06  &  2387.13  &  1.47~e+06  &  3463.57  &  \textbf{1.42~e+06}  &  2390.27  & 2.15~e+06 & 14400  \\
        20 & 0.5 & 2 &  \textbf{0.71~e+06}  &  1931.70  &  0.73~e+06  &  2503.33  &  0.75~e+06  &  1813.50  & 9.50~e+05 & 14400  \\
        20 & 0.5 & 3 &  2.18~e+05  &  1625.86  &  \textbf{2.09~e+05}  &  2383.20  &  2.15~e+05  &  1406.08  & 3.12~e+05 & 14400  \\
        20 & 0.5 & 4 &  \textbf{7.39~e+05}  &  2788.98  &  7.71~e+05  &  2917.41  &  8.37~e+05  &  2186.68  & 8.35~e+05 & 14400  \\
        20 & 0.5 & 5 &  0.96~e+06  &  2011.54  &  0.97~e+06  &  2678.38  &  \textbf{0.90~e+06}  &  1676.48  & 1.07~e+06 & 14400  \\
        20 & 0.5 & 6 &  3.24~e+06  &  2536.07  &  2.89~e+06  &  3327.69  &  \textbf{2.65~e+06}  &  2685.08  & 3.95~e+06 & 14400  \\
        20 & 0.5 & 7 &  2.08~e+05  &  2199.30  &  \textbf{2.01~e+05}  &  3456.32  &  2.03~e+05  &  2229.30  & 2.56~e+05 &  14400  \\
        20 & 0.5 & 8 &  2.43~e+05  &  1641.64  &  2.42~e+05  &  2354.32  &  \textbf{2.41~e+05}  &  1514.79  & 3.37~e+05 &  14400  \\
        20 & 0.5 & 9 &  6.20~e+05  &  2529.35  &  6.24~e+05  &  3374.94  &  5.80~e+05  &  2948.75  & \textbf{5.51~e+05} & 14400 \\
    \hline
        20 & 0.75 & 0 &  \textbf{6.90~e+04}  &  1783.51  &  7.30~e+04  &  3277.64  &  7.30~e+04  &  1188.81  & 7.86~e+04 &  14400\\
        20 & 0.75 & 1 &  1.12~e+05  &  1900.06  &  1.12~e+05  &  3649.95  &  \textbf{1.12~e+05}  &  1753.16  & 1.18~e+05 &  14400\\
        20 & 0.75 & 2 &  \textbf{0.83~e+05}  &  1691.69  &  0.87~e+05  &  2724.63  &  0.88~e+05  &  1159.22  & 1.01~e+05 &  14400\\
        20 & 0.75 & 3 &  \textbf{7.00~e+04}  &  2235.44  &  7.50~e+04  &  2782.49  &  7.90~e+04  &  1148.20  & 9.35~e+04 &  14400\\
        20 & 0.75 & 4 &  \textbf{0.99~e+05}  &  1735.64  &  1.00~e+05  &  2965.26  &  1.02~e+05  &  1159.36  & 1.06~e+05 &  14400\\
        20 & 0.75 & 5 &  \textbf{0.97~e+05}  &  1963.95  &  0.97~e+05  &  3402.85  &  0.97~e+05  &  2020.67  & 1.01~e+05 &  14400\\
        20 & 0.75 & 6 &  \textbf{0.91~e+05}  &  1508.97  &  0.93~e+05  &  3555.80  &  0.92~e+05  &  1476.82  & 9.74~e+04 &  14400\\
        20 & 0.75 & 7 &  \textbf{0.89~e+05}  &  1710.47  &  0.89~e+05  &  3122.01  &  0.90~e+05  &  1367.58  & 1.01~e+05 &  14400\\
        %20 & 0.75 & 8 &  0.94~e+05  &  516.30  &   &   &  0.97~e+05  &  1427.10  & 1.09~e+05 & 14400 \\
        20 & 0.75 & 8 &  \textbf{0.98~e+05}  &  1776.35  &  0.99~e+05  &  3183.02  &  1.01~e+05  &  2031.30  & 1.11~e+05 &  14400\\
        20 & 0.75 & 9 &  \textbf{0.99~e+05}  &  1779.24  &  1.02~e+05  &  4046.86  &  1.02~e+05  &  1470.36  & 1.05~e+05 & 14400 \\
    \hline
        30 & 0.5 & 0 &  1.97~e+05  &  8131.35  &  1.97~e+05  &  9303.11  &  \textbf{1.84~e+05}  &  6661.51  & 1.98~e+05 & 14400 \\
        30 & 0.5 & 1 &  2.09~e+05  &  6741.38  &  \textbf{2.05~e+05}  &  8337.86  &  2.06~e+05  &  6223.99  & 2.36~e+05 & 14400 \\
        30 & 0.5 & 2 &  1.62~e+05  &  6310.00  &  1.63~e+05  &  5466.08  &  \textbf{1.61~e+05}  &  2403.93  & 1.80~e+05 & 14400 \\
        30 & 0.5 & 3 &  \textbf{1.76~e+05}  &  7258.44  &  1.77~e+05  &  8490.38  &  1.83~e+05  &  5364.12  & 2.38~e+05 & 14400 \\
        30 & 0.5 & 4 &  3.04~e+05  &  8005.14  &  3.22~e+05  &  8665.91  &  \textbf{2.89~e+05}  &  6750.19  & 3.49~e+05 & 14400 \\
        30 & 0.5 & 5 &  \textbf{1.88~e+05}  &  7435.17  &  1.88~e+05  &  7963.70  &  1.91~e+05  &  4415.26  & 2.26~e+05 & 14400 \\
        30 & 0.5 & 6 &  1.88~e+05  &  7244.06  &  1.95~e+05  &  8679.78  &  \textbf{1.88~e+05}  &  5036.13  & 2.33~e+05 & 14400 \\
        30 & 0.5 & 7 &  \textbf{3.82~e+05}  &  8138.68  &  4.05~e+05  &  9202.13  &  4.20~e+05  &  460.71  & 4.39~e+05 & 14400 \\
        30 & 0.5 & 8 &  \textbf{1.51~e+05}  &  5786.28  &  1.52~e+05  &  7515.79  &  1.57~e+05  &  4728.17  & 1.86~e+05 & 14400 \\
        30 & 0.5 & 9 &  2.78~e+05  &  7180.90  &  2.73~e+05  &  7172.85  &  \textbf{2.59~e+05}  &  4259.48  & 3.30~e+05 & 14400 \\
    \hline
        30 & 0.75 & 0 &  1.04~e+05  &  5948.39  &  \textbf{1.02~e+05}  &  11937.47  &  1.05~e+05  &  3976.56  & 1.15~e+05 &  14400 \\
        30 & 0.75 & 1 &  1.23~e+05  &  5546.67  &  \textbf{1.22~e+05}  &  10765.75  &  1.28~e+05  &  3371.02  & 1.37~e+05 &  14400 \\ 
        30 & 0.75 & 2 &  0.98~e+05  &  7392.68  &  \textbf{0.96~e+05}  &  12556.06  &  1.01~e+05  &  3642.89  & 1.14~e+05 & 14400 \\
        30 & 0.75 & 3 &  1.03~e+05  &  6083.00  &  \textbf{0.99~e+05}  &  12773.28  &  1.04~e+05  &  2800.32  & 1.08~e+05 & 14400 \\
        30 & 0.75 & 4 &  0.98~e+05  &  6049.69  &  \textbf{0.95~e+05}  &  12036.30  &  1.02~e+05  &  3160.86  & 1.15~e+05 & 14400 \\
        30 & 0.75 & 5 &  1.04~e+05  &  6026.44  &  \textbf{1.02~e+05}  &  13553.03  &  1.06~e+05  &  3695.71  & 1.13~e+05 & 14400 \\
        30 & 0.75 & 6 &  0.99~e+05  &  6660.11  &  \textbf{0.98~e+05}  &  11316.65  &  0.99~e+05  &  3509.40  & 1.11~e+05 & 14400 \\
        30 & 0.75 & 7 &  1.07~e+05  &  6519.37  &  \textbf{1.05~e+05}  &  11515.46  &  1.09~e+05  &  3668.31  & 1.27~e+05 & 14400 \\
        30 & 0.75 & 8 &  \textbf{0.93~e+05}  &  6524.25  &  0.94~e+05  &  11144.66  &  0.96~e+05  &  3562.64  & 1.06~e+05 & 14400 \\
        30 & 0.75 & 9 &  0.89~e+05  &  7337.74  &  \textbf{0.88~e+05}  &  12763.55  &  0.93~e+05  &  4114.14  & 1.04~e+05 & 14400 \\
    \hline
    \end{tabular}
    \end{adjustwidth}
    \caption{Results on random instances with n nodes randomly generated in U(0,1), p probability of an edge with an additional factor based on distance, $i$ instance number. An additional $\frac{n}{10}$ nodes are marked as exits on the boundary of the unit square. These also serve as FR entry points. Capacities of edges are N(50,20). All edges have two lanes in both directions. The population at each node to exit follows a N(100,10). For each instance, the entry for the method with the best performance is in bold text.}
    \label{tab:random_instances_comparison}
\end{table}

\subsubsection{Varying the parameter $n_{paths}$}

To understand the impact of increasing the number of Graver elements, we examine a single random instance ($n=30$, $p=0.75$, $i=1$) in further detail. This instance is chosen because the solution from GAGA-Yen's algorithm is slightly higher when compared to the GAGA-SA-10000 result. In particular, we explore the effect of changing $n_{paths}$ and compare GAGA-SA-10000 and GAGA-Yens results. We vary the $n_{paths}$ from $20$ to $40$ in steps of $5$. We show the solution and time to solution in Table~\ref{tab:random_instance_graver_variation}. The solution at the end of the GAGA run (GAGA-\textit{only}) is shown in brackets. The table shows that the GAGA-\textit{only} result is comparable with the GAGA-\textit{LeBlanc} for the GAGA-SA-10000. The effects of changing $n_{paths}$ are negligible. On the other hand, if we examine the GAGA-Yens results, we can see that the GAGA-\textit{only} results are significantly larger than the GAGA-\textit{LeBlanc} values. While the effect of $n_{paths}$ is not significant on the final solution, we can see that the GAGA-\textit{only} solution decreases with increasing $n_{paths}$ (which also leads to an increase in the number of Graver elements, which we denote as $n_{graver}$).

\begin{table}[ht]
\fontsize{9}{9}\selectfont
\renewcommand{\arraystretch}{1.2}
\centering
%\begin{adjustwidth}{-2cm}{}
    \begin{tabular}{|c|c|c|c|c|c|c|c|c|c|c|c|}
    \hline
        \multicolumn{1}{|c|}{Instance} & \multicolumn{3}{|c|}{GAGA-SA-10000} & \multicolumn{3}{|c|}{GAGA-Yens} & \multicolumn{2}{|c|}{Branch-And-Bound} \\
    \hline
        $n_{paths}$ & $n_{graver}$ & obj & time (s) & $n_{graver}$ & obj & time (s)& obj & time (s) \\
    \hline

        20 &  10096  &  1.23~e+05  &  8850.87  &  7432  &  1.30~e+05  &  1662.32  & 1.41~e+05 &  14400 \\
         &    &  (1.24~e+05)  &  (7224.89 + 138.48  &    &  (3.99~e+05)  &  (201.57 + 127.22 & &  \\
         &    &    &  + 1321.51 + 165.99)  &    &    &   + 1176.85 + 156.67)  & &  \\
 
        25 &  15854  &  1.22~e+05  &  10763.01  &  11484  &  1.25~e+05  &  3727.72  & &  \\
         &    &  (1.24~e+05)  &  (7224.89 + 833.29  &    &  (2.81~e+05)  &  (201.57 + 584.36   & &  \\
        &    &    &  + 2442.32 + 161.50)  &    &    &  + 2781.59 + 160.19)  & &  \\
 
        30 &  22942  &  1.22~e+05  &  14547.90  &  16187  &  1.28~e+05  &  6534.44  & &  \\
         &    &  (1.23~e+05)  &  (7224.89 + 3058.06   &    &   (2.67~e+05) &  (201.57 + 1907.81   & &  \\
          &    &    &   + 4102.57 + 162.37)  &    &    &   + 4260.08 + 164.97)  & &  \\
  
        35 &  31382  &  1.23~e+05  &  22431.84  &  22019  &  1.28~e+05  &  14720.39  & &  \\
         &    & (1.24~e+05)   &  (7224.89 + 8714.43   &    &   (1.91~e+05) &  (201.57 + 5160.53   & &  \\
         &    &    &   + 6333.83 + 158.67)  &    &    &   + 9764.47 + 165.97)  & &  \\
 
        40 &  41065  &  1.22~e+05  &  38294.72  &  27886  &  1.29~e+05  &  26939.96  & &  \\
         &    &  (1.23~e+05)  &  (7224.89 + 20923.10   &    & (1.82~e+05)   &  (201.57 + 11832.86   & &  \\
          &    &    &   + 9986.53+160.18)  &    &    &   + 14740.83+164.69)  & &  \\
\hline
    \end{tabular}
    %\end{adjustwidth}
    \caption{For a single graph instance ($n=30, p=0.75, i=1$), comparison between simulated annealing with 10000 samples (GAGA-SA-10000) and Yens k-shortest path method (GAGA-Yens) with a varying number of paths ($k=20, 25, 30, 35, 40$). The quantity in brackets of the column labeled, \textit{"obj"} of GAGA-SA-10000 and GAGA-Yens refers to the solution at the end of the GAGA run only (i.e., GAGA-\textit{only} solution). The time taken for each method of the GAGA algorithm (feasible solutions calculation, Graver basis computation, Graver walk, and \textit{LeBlanc} algorithm) are respectively shown in the brackets of the "time(s)" column.}
    \label{tab:random_instance_graver_variation}
\end{table}

\subsubsection{Experimenting with interchanging initial seeds and Graver basis from different path generation methods}

To further understand the effects of the initial seeds and partial Graver basis elements, we run the GAGA algorithm using initial seeds generated from GAGA-SA-10000/GAGA-Yens and the partial Graver basis from GAGA-Yens/GAGA-SA-10000. The results are compared in Table~\ref{tab:random_instance_graver_variation_same_initial}. The final GAGA-\textit{LeBlanc} values are comparable in all cases. Instead, if we compare the GAGA-\textit{only} solutions, we can see that the best result is obtained when using GAGA-SA-10000 solutions both for generating the initial solutions and the Graver elements. If we use initial seeds from GAGA-SA-10000 solutions but generate the Graver basis elements using GAGA-Yens, the GAGA-\textit{only} results are worse by 2 orders of magnitude. This can be attributed to the Graver basis elements from GAGA-Yens not being able to cover the domain of the paths found in the GAGA-SA-10000 initial feasible solutions.

\begin{table}
\centering
    \begin{tabular}{|c|c|c|c|c|}
    \hline
        \multicolumn{3}{|c|}{Instance} & \multicolumn{1}{|c|}{GAGA-SA-10000/GAGA-Yens} \\
    \hline
        Initial seed Generator & Graver generator & $n_{graver}$ & obj\\
    \hline
       SA  &  SA  & 15854  &  1.23~e+05 (1.23~e+05) \\
      Yens &  SA  & 15854  &  1.25~e+05 (1.36~e+05) \\
      Yens & Yens & 22019  &  1.28~e+05 (1.91~e+05) \\
       SA  & Yens & 22019  &  1.25~e+05 (206.17~e+05) \\
    \hline
    \end{tabular}
    \caption{For a single graph instance ($n=30, p=0.75, i=1$), comparison between simulated annealing (GAGA-SA-10000) and Yen's k-shortest paths (GAGA-Yens) with the same initial conditions and Graver basis interchanged. The quantity in brackets of the column, \textit{"obj"} refers to the solution at the end of the GAGA run only (GAGA-\textit{only} solution).}
    \label{tab:random_instance_graver_variation_same_initial}
\end{table}

Given that the quality of solutions is consistently better with the GAGA-SA-10000 feasible solutions in $n=30$ and $p=0.75$ experiments, and we had observed that GAGA-{\emph only} solutions are better for  $n=10000$ (versus $n=1000$ for SA),  in Section \ref{sec: Turkish instances}, we decided to use feasible paths generated by simulated annealing run with $n_{samples}=10000$.

\subsection{Case study: Istanbul instances} \label{sec: Turkish instances}

We generate three sample instances based on the city of Istanbul, Türkiye. This enables us to assess GAGA's performance on a real-life urban road network. The case study instance generation and the performance comparison of GAGA and the branch and bound algorithm are detailed in Sections \ref{sec: InstanveGen} and \ref{sec: GAGAvsBB}.

\subsubsection{Instance generation} \label{sec: InstanveGen}
To demonstrate the impact of FR lane reservation decisions on evacuees' travel times, we designate a realistic case study on a specific region in Istanbul, including the Avcilar district and a neighborhood around the Küçükçekmece Lake. Istanbul is subject to major earthquake risk due to the fault line passing under the sea to the city's south \cite{erdik2008}, and Avcilar is one of the most risky areas. The region of study is home to nearly 500,000 residents \cite{TUIK2023}. 

\begin{figure}[ht]
    \centering
    \includegraphics[width=0.5\linewidth]{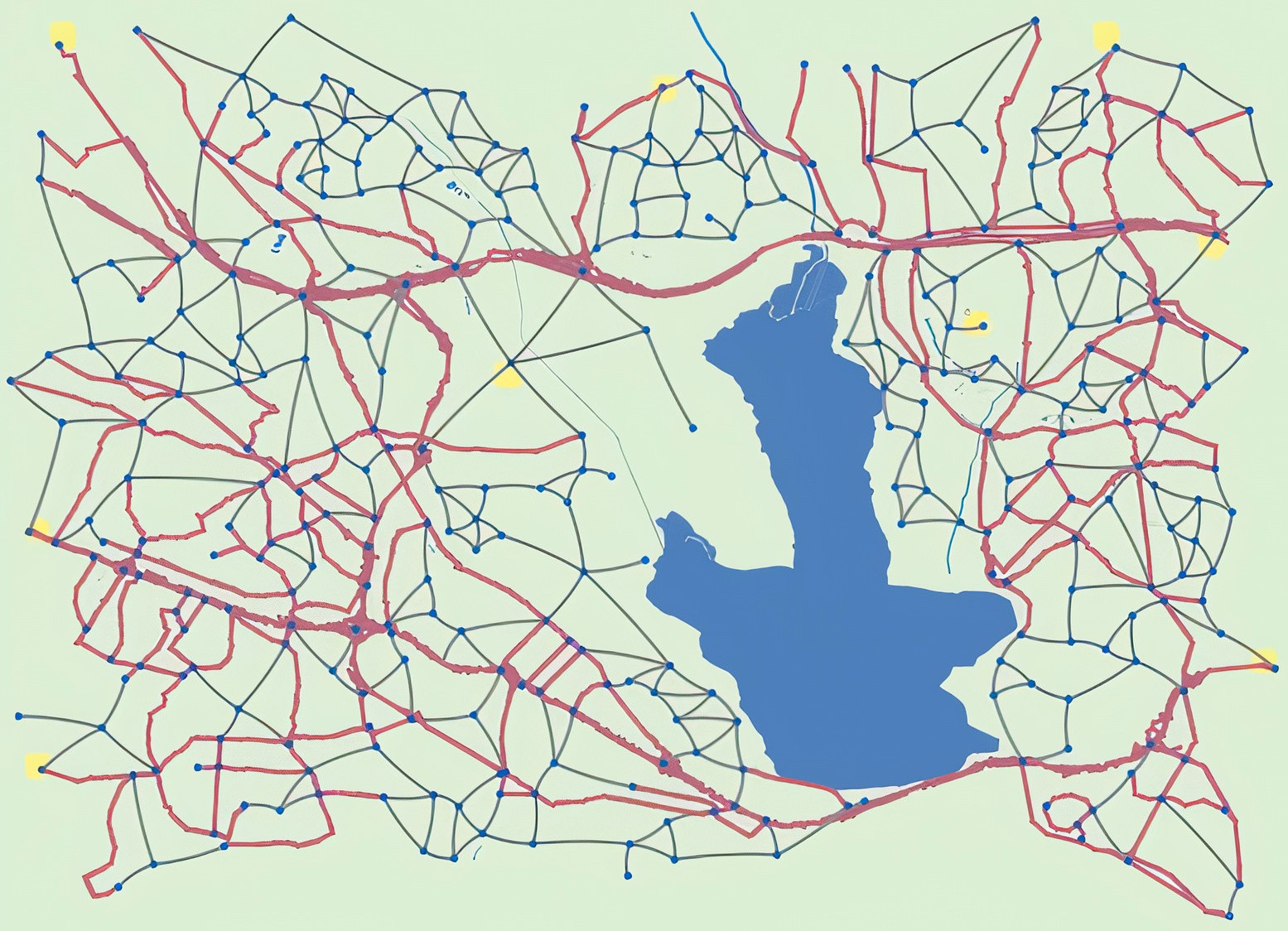}
    \caption{The network associated with our case study instances. FR lanes can be reserved along the thick pink edges. The yellow circles stand for exit nodes.}
    \label{fig:region}
\end{figure}

To generate the region's road network, we leveraged Istanbul's OpenStreetMap (OSM) dataset, focusing on the road infrastructure information layer. This layer is analyzed in ArcGIS 10.5, where links are systematically classified into two groups based on their road-type features. Among the 15 distinct features for road records, the first class encompasses links with road-type features such as \textit{primary}, \textit{secondary}, \textit{motorway}, and \textit{trunk}. These links are precisely mapped in the network and hold the potential for hosting an FR lane reservation. On the other hand, the second class includes lanes with \textit{tertiary}, \textit{residential}, \textit{service}, \textit{track}, and \textit{living street} road-type features. While these links cannot accommodate FR lane reservations due to capacity and lane limitations, they remain accessible for evacuees and FRs to navigate the network. Given the substantial number of links falling into the second class, the selected arcs comprise only main roads featuring two or more lanes on both sides. Finally, link types such as \textit{cycleway}, \textit{pedestrian}, \textit{footway}, \textit{steps}, and \textit{path} are excluded from the network, as they are impassable for both evacuees and FRs. The resulting network, depicted in Figure \ref{fig:region}, consists of 179 vertices, each representing a residential area or an origin node for the FRs, and 234 arcs representing the road infrastructure connecting the vertices in both directions. The thick pink and the thin grey links correspond to the first and second categories of edges, respectively. The yellow patches represent exit nodes. An exit node situated within the region indicates a safe area.

\cite{yabe2019cross} studied a cross-comparative analysis of individuals' evacuation behaviors in a post-earthquake situation. The authors investigated the GPS trajectories of more than one million anonymized mobile phone users whose positions were tracked for two months before and after four of the major earthquakes that occurred in Japan. They argued that in the face of an earthquake with a 7 seismic intensity, nearly 60\% of individuals strive to evacuate the region immediately. This result is, however, regardless of the earthquake hitting time, which can affect the magnitude of the evacuation. Therefore, jointly considering such an evacuation ratio with uncertainties in evacuees' travel volume, we generated three instances on the same network by varying the evacuation demand level. Based on the population size $p_i$ at each node $i$ provided in the \cite{Avcilar2023} report, 4-person vehicle capacities, and the evacuation ratio reported by \cite{yabe2019cross}, we assume among the three instance sets 1, 2, and 3, the number of evacuating vehicles $d_i$ emanating from each demand node $i$ is measured by $\big|0.25 p_i \times U[0.1,0.3]\big|$, $\big|0.25 p_i \times U[0.3,0.5]\big|$, and $\big|0.25 p_i \times U[0.5,0.7]\big|$, respectively. These uniform distributions $U$ account for demand randomness and capture different evacuation ratios depending on the occurrence time of an earthquake with a seismic intensity of 7.

\subsubsection{Comparing GAGA and BB} \label{sec: GAGAvsBB}

Within the designated network in Section \ref{sec: Turkish instances}, there exists $154$ nodes with evacuation demand and $69$ nodes with FR demands, varying in demand levels. We generate feasible paths from each demand node to the exit nodes, taking into account the traffic equilibrium. Subsequently, we construct the partial Graver basis for both FR and evacuee paths. We compare the GAGA results using the `GAGA-SA-10000' results with the branch-and-bound results. In the GAGA algorithm, four different settings are explored, as discussed below.

\begin{itemize}
    \item \textit{Unnormalized}: The evacuee demand is fixed as it is. The Graver walk in the inner user-equilibrium problem is terminated only when no further augmentation is possible at that solution.
    \item \textit{Unnormalized, Tol = 1e-3}: The evacuee demand is fixed as it is. The Graver walk in the inner problem is terminated when the fractional difference in optima before and after one full pass over all the Graver elements is less than $10^{-3}$.
    \item \textit{Normalized}: The evacuee demand is normalized such that the maximum demand at any node is 100. The optimization in the inner problem is terminated only when no further augmentation is possible at that solution.
    \item \textit{Normalized, Tol=1e-3}: The tolerance before and after one full pass over all the Graver elements is less than $10^{-3}$.
\end{itemize}

The Branch-and-Bound (BB) algorithm is run up to $72$ hours. The Graver walk step in the GAGA algorithm is also given a time limit of $72$ hours. Except for the \textit{Unnormalized} setting, the total run time of GAGA for the remaining $3$ settings is less than that of BB for all the $3$ instances. For the \textit{Unnormalized} setting, the total run time for the $1^{st}$ instance is less than that of BB. In the remaining $2$ instance, we find that the Graver augmentation time increases due to increased evacuee demand. Hence, the Graver walk is not fully completed within the specified time limit for all the initial seeds. We note that all the runs for these instances have been performed on Processor-G. 

\begin{figure*}
\centering
     \includegraphics[width=0.7\textwidth]{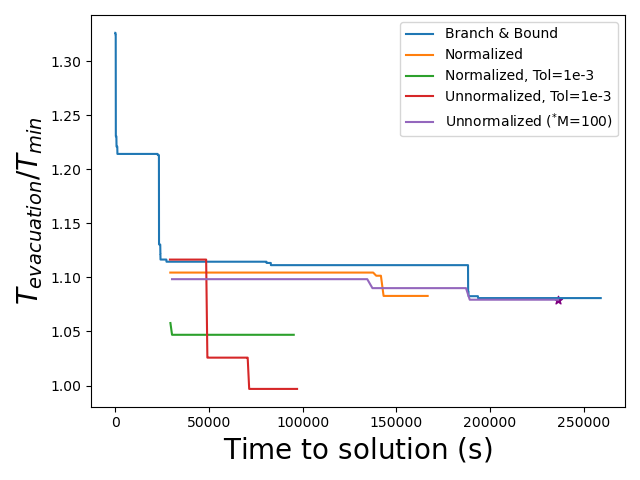}
     \includegraphics[width=0.7\textwidth]{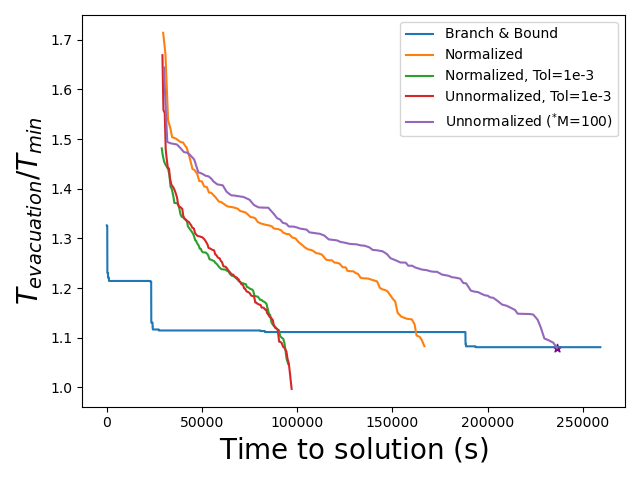}
     \caption{The evacuation times as a function of time to solution for Turkish Graph Instance-1 computed using Graver walk with varying settings (on evacuee demand normalization and tolerance threshold in inner user equilibrium problem) compared against Branch and Bound solution. For each run, the final FR paths (local optimal solutions for each seed) obtained after the Graver walk are all found to be unique. On average, the inner user-equilibrium level with tolerance threshold is faster by a factor of $\sim 3.85$ than without for the \textit{Unnormalized} setting and $\sim 2.5$ for a \textit{Normalized} setting. \textit{Top panel:} The best solution at the given time with the order of the seeds is the same as in the computational experiment. \textit{Bottom panel:} The order of seeds is based on the decreasing order of the objective function.}
     \label{fig:T_instance1}
\end{figure*}

\begin{figure*}
    \centering
     \includegraphics[width=0.7\textwidth]{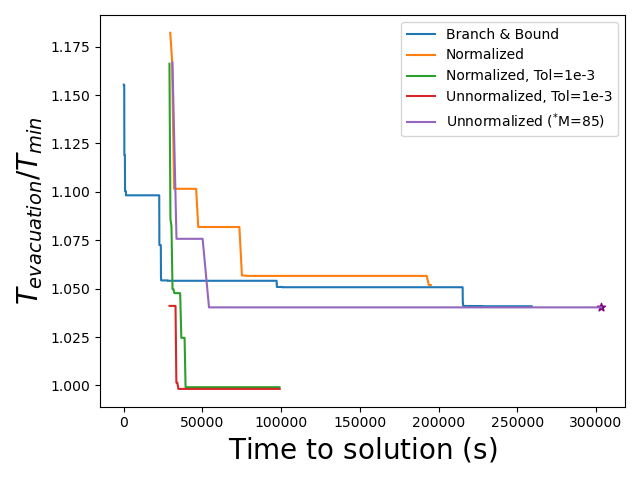}
     \includegraphics[width=0.7\textwidth]{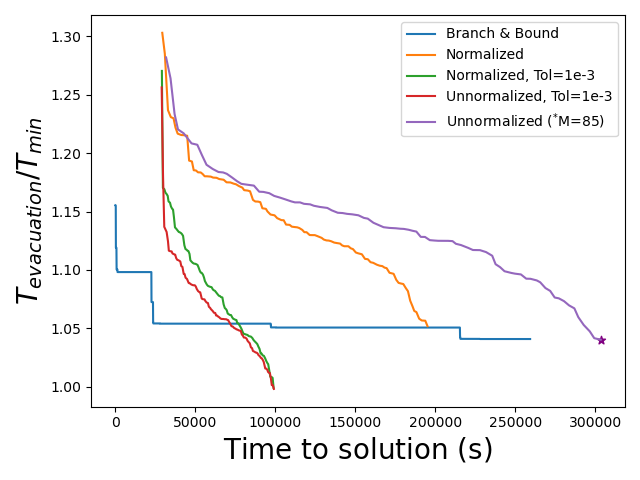}
     \caption{The evacuation times as a function of time to solution for Turkish Graph Instance-2 computed using Graver walk with varying settings (on demand normalization and tolerance) compared against Branch and Bound solution. For the GAGA run with \textit{Unnormalized}, only (*M=85) seeds are used due to the time limit. The remaining GAGA runs use 100 seeds. For each run, the final FR paths (local optimal solutions for each seed) obtained after the Graver walk are all found to be unique. On average, the inner user-equilibrium level with tolerance threshold is faster by a factor of $\sim 6.1$ than without for the \textit{Unnormalized} setting and $\sim 3$ for a \textit{Normalized} setting.}
     \label{fig:T_instance2}
\end{figure*}

\begin{figure*}
    \centering
     \includegraphics[width=0.7\textwidth]{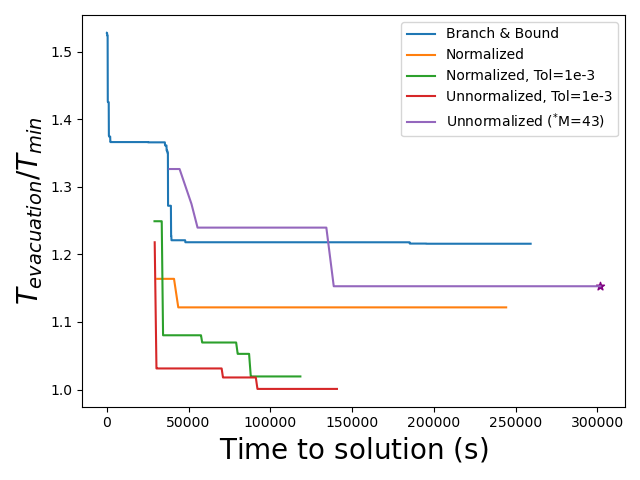}
     \includegraphics[width=0.7\textwidth]{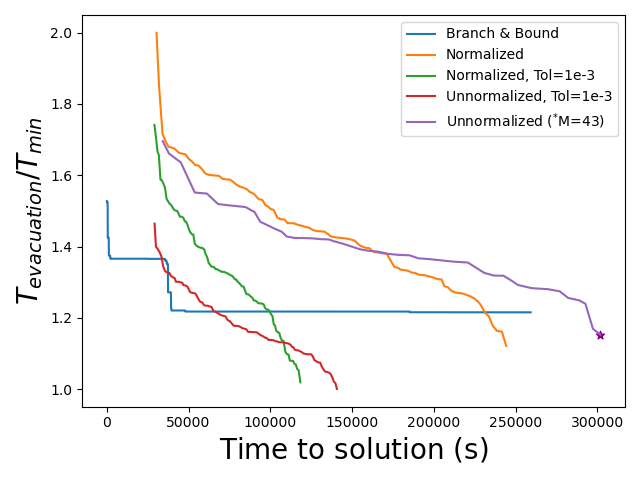}
     \caption{The evacuation times as a function of time to solution for Turkish Graph Instance-3 computed using Graver walk with varying settings (on demand normalization and tolerance) compared against Branch and Bound solution. For the GAGA run with \textit{Unnormalized} setting, only (*M=43) seeds are used due to the time limit. The remaining GAGA runs use 100 seeds. For each run, the final FR paths (local optimal solutions for each seed) obtained after the Graver walk are all found to be unique. On average, the inner user-equilibrium level with tolerance threshold is faster by a factor of $\sim 9.56$ than without for the \textit{Unnormalized} setting and $\sim 3.44$ for a \textit{Normalized} setting.}
     \label{fig:T_instance3}
\end{figure*}

The results comparing each of the $4$ settings and comparing with Branch-and-Bound are tabulated in Tables ~\ref{tab:t_graph2}, ~\ref{tab:t_graph_normalized}, ~\ref{tab:t_graph_tole}, ~\ref{tab:t_graph_tole_norm}. The run times for each of the three instances are plotted in Figures~\ref{fig:T_instance1}, ~\ref{fig:T_instance2}, ~\ref{fig:T_instance3}. In all the figures, there is a pre-Graver-walk time for each of the $4$ settings to compute the feasible solutions and partial Graver basis. The BB run instead starts from $t=0 s$. 

In each figure, the top panel shows the time to solution and the optimal solution at that time, following the order of initial seeds generated randomly from the Simulated annealing solutions. The bottom panel shows the same results but follows the ``worst case" ordering where the seed with the worst optimal solution comes first, and the optimal solution decreases as we proceed to the next seeds.

In Figure~\ref{fig:T_instance1}, we show the results for the first instance. We can see that the optimal solution from Branch and bound is similar to the GAGA result for the \textit{Unnormalized} and \textit{Normalized} settings, but larger than the results when including a tolerance threshold setting by $\sim 10\%$. Comparing amongst the GAGA runs, we notice that the time to solution is much longer for the \textit{Unnormalized} and \textit{Normalized} setting with no tolerance setting. The best solution is obtained in the \textit{Unnormalized} setting with a tolerance threshold setting included. This setting requires much less time and is also comparable to the Normalized setting with a tolerance threshold. In Figure~\ref{fig:T_instance2} and Figure~\ref{fig:T_instance3}, we show the same results for the $2^{nd}$ and $3^{rd}$ instances. We again find that the best GAGA results are obtained with the \textit{Unnormalized} and \textit{Normalized} settings with a tolerance threshold and are better than BB by $5\%$ and $20 \%$ respectively for the two instances. The BB results are similar to the GAGA runs with only the \textit{Unnormalized} and \textit{Normalized} settings. So, we can conclude that running the GAGA algorithm in an \textit{Unnormalized} or \textit{Normalized} and including a tolerance criterion in the inner user equilibrium problem tends to produce better quality solutions when compared to BB and in a shorter time interval.

\begin{table}
\centering
    \begin{tabular}{|c|c|c|c|c|c|c|}
    \hline
        \multicolumn{3}{|c|}{Instance}  & \multicolumn{2}{|c|}{GAGA-SA-10000} & \multicolumn{2}{|c|}{Branch-And-Bound} \\
    \hline
        n & m & l & obj & time (s)  & obj & time (s) \\ \hline
    \hline
        179 & 234 & 1 &  1.54~e+06 (M=100) &  236438.77  &   1.55~e+06 &  259200 \\
        %179 & 234 & 1 &  1.54~e+06 (M=100) &  236438.77  &   1.59~e+06 &  92000 \\        
        % &  &  &    &    &  &  (23482.36 + 4847.92 +  84542.02) &    &   &   &  \\
        179 & 234 & 2  & 2.56~e+07 (M=85) &  303312.19 &  2.56~e+07 &  259200 \\
        %179 & 234 & 2  & 2.56~e+07 (M=85) &  303312.19 &  2.59~e+07 &  92000 \\
         %&  &  &    &    &  &  (23482.36 + 4847.92 + 89299.19) &    &   &   &  \\
        179 & 234 & 3  & 1.16~e+08 (M=43) &  301813.49 &   1.23~e+08 &  259200 \\
        %179 & 234 & 3  & 1.16~e+08 (M=43) &  301813.49 &   1.23~e+08 &  92000 \\
         %&  &  &    &    &  &  (23482.36 + 4847.92 + 83697.59) &    &   &   &  \\
    \hline
    \end{tabular}
    \caption{Comparison of the solution output by GAGA and time to solution with Branch-And-Bound on 3 Turkish graph instances (Avcilar graph with different demands). The GAGA algorithm is run in the \textit{Unnormalized} setting without any normalization of evacuee demand, and there is no adjustment of tolerance in the inner problem.}
    \label{tab:t_graph2}
\end{table}

%\begin{figure*} 
%     \includegraphics[width=1.0\textwidth]{figures/fig_Treal_soltype_1.png}
%     \caption{The evacuation times for different instances with varying demands plotted against the time to solution. Here, the criterion for graver walk termination in the inner level is fixed at a tolerance of $10^{-3}$.}
%     \label{fig:T_tol}
%\end{figure*}

\begin{table}
\centering
    \begin{tabular}{|c|c|c|c|c|c|c|}
    \hline
        \multicolumn{3}{|c|}{Instance}  & \multicolumn{2}{|c|}{GAGA-SA-10000}  & \multicolumn{2}{|c|}{Branch-And-Bound} \\
    \hline
        n & m & l & obj & time (s)  & obj & time (s)  \\ \hline
    \hline
        179 & 234 & 1 &  1.55~e+06 (M=100) &  166784.54  &  1.55~e+06 &  259200 \\
        179 & 234 & 2  & 2.59~e+07 (M=100) &  195092.25  &  2.56~e+07 &  259200 \\
        179 & 234 & 3  & 1.13~e+08 (M=100) &  244233.26  &  1.23~e+08 &  259200 \\
        %179 & 234 & 1 &  1.55~e+06 (M=100) &  166784.54  &  1.59~e+06 &  92000 \\
        %179 & 234 & 2  & 2.59~e+07 (M=100) &  195092.25  &  2.59~e+07 &  92000 \\
        %179 & 234 & 3  & 1.13~e+08 (M=100) &  244233.26  &  1.23~e+08 &  92000 \\
    \hline
    \end{tabular}
    \caption{Comparison of the solution output by GAGA and time to solution with Branch-And-Bound on 3 Turkish graph instances (same as in Table~\ref{tab:t_graph2}) with the GAGA run with \textit{Normalized} setting. The maximum evacuee demand in each instance is normalized to 100, and an equivalent factor reduces capacity before Graver walks. There is no adjustment to tolerance in the inner problem.}
    \label{tab:t_graph_normalized}
\end{table}

\begin{table}
\centering
    \begin{tabular}{|c|c|c|c|c|c|c|}
    \hline
        \multicolumn{3}{|c|}{Instance}  & \multicolumn{2}{|c|}{GAGA-SA-10000}  & \multicolumn{2}{|c|}{Branch-And-Bound} \\
    \hline
        n & m & l & obj & time (s)  & obj & time (s)  \\ \hline
    \hline
        179 & 234 & 1  &   1.43~e+06 (M=100) &  97055.87  &  1.55~e+06 &  259200 \\
        179 & 234 & 2  &  2.46~e+07  (M=100) &  99077.13 &  2.56~e+07 &  259200 \\
        179 & 234 & 3  &  1.01~e+08  (M=100) &  140761.56 &  1.23~e+08 &  259200 \\
    \hline
    \end{tabular}
    \caption{Comparison of the solution output by GAGA and time to solution with Branch-And-Bound on 3 Turkish graph instances (same as in Table~\ref{tab:t_graph2}) with the GAGA run with \textit{Unnormalized, Tol = 1e-3} setting. There is no normalization of evacuee demand, and the criterion for Graver walk termination in the inner user equilibrium level is a tolerance threshold of $10^{-3}$.}
    \label{tab:t_graph_tole}
\end{table}

\begin{table}
\centering
    \begin{tabular}{|c|c|c|c|c|c|c|}
    \hline
        \multicolumn{3}{|c|}{Instance}  & \multicolumn{2}{|c|}{GAGA-SA-10000}  & \multicolumn{2}{|c|}{Branch-And-Bound} \\
    \hline
        n & m & l & obj & time (s)  & obj & time (s)  \\ \hline
    \hline
        179 & 234 & 1  &   1.50~e+06 (M=100) &  95280.47  &  1.55~e+06 &  259200 \\
        179 & 234 & 2  &  2.46~e+07  (M=100) &  98936.04 &  2.56~e+07 &  259200 \\
        179 & 234 & 3  &  1.03~e+08  (M=100) &  118397.11 &  1.23~e+08 &  259200 \\
    \hline
    \end{tabular}
    \caption{Comparison of the solution output by GAGA and time to solution with Branch-And-Bound on 3 Turkish graph instances (same as in Table~\ref{tab:t_graph2}) with the GAGA run with \textit{Normalized, Tol = 1e-3} setting. Here, the maximum evacuee demand in each instance is normalized to 100, and an equivalent factor reduces capacity before Graver walks. The criterion for Graver walk termination in the inner problem is a tolerance threshold of $10^{-3}$.}
    \label{tab:t_graph_tole_norm}
\end{table}

\section{Conclusion}

Following a disaster, the immediate priority lies in locating and rescuing individuals who may be trapped or injured while ensuring the other entities can travel to safe locations on their own. This is further effective in case the available and restricted lanes are properly labeled, having been pre-determined and pre-communicated. However, identifying such roads requires a careful assessment of the tradeoff between expediting the FR responses and reducing the travel time for the disaster-affected people. This preparedness activity is parallel to other efforts, including establishing medical facilities, providing emergency medical care, and ensuring the availability of crucial medical supplies. 

%Deploying law enforcement personnel to affected areas is vital to ensure the safety of both responders and survivors. These critical activities are contingent on accessible roads, which inevitably face heightened traffic congestion due to the mobile public utilizing the same road network system. To facilitate faster access for the FRs, some road segments can be specifically reserved for their use. Identifying these emergency routes before a disaster has the advantage of informing the public beforehand, acting faster after the disaster, and strengthening them structurally against potential disasters in the preparedness stage. 

We define the first responder network design problem as a bilevel optimization problem, where the arcs reserved for the first responders are selected at the first level to minimize total travel time, which depends on the second level problem where the public travels between their origin-destination pairs over the eligible roads, adopting a selfish routing principle. We propose a novel solution methodology for this difficult bilevel nonlinear network design problem. Unlike the existing methodologies for similar network design problems, we apply a quantum-inspired algorithm. The algorithm takes advantage of the Graver Augmented Multi-Seed Algorithm (GAMA) and follows a bi-level nested GAMA within the GAMA algorithm that we name GAGA. We test several variations of GAGA on instances related to a set of random graph instances and three scenarios on a predicted Istanbul earthquake. In most instances, we elicit superior solution quality and run-time performance of GAGA compared to a state-of-the-art exact algorithm for a traditional formulation.

We hope that our promising results encourage applying quantum-inspired methods to other disaster preparedness problems with combinatorial aspects and computationally expensive nonlinearity in the objective function and also explore further development of this novel methodology to such complex models in other application domains. %The results of our study not only point to promising new research directions to extend the use of quantum and quantum-inspired algorithms to solve complex decision problems but also inform hardware developers to identify the key areas of development for enhancing the capabilities of quantum computing technologies for optimization. 

\vspace{0.1in}

\noindent{\bf Acknowledgements.} ST and AT would like to acknowledge Raytheon BBN for their support through a CMU-BBN contract as part of a DARPA  project on Quantum Inspired Classical Computing (QuICC). For AK, this material is also based upon work supported by the National Science Foundation Graduate Research Fellowship Program under Grant No. DGE1745016, DGE2140739. Any opinions, findings, conclusions, or recommendations expressed in this material are those of the author(s) and do not necessarily reflect the views of the National Science Foundation.

\vspace{0.1in}

\section*{Appendix A: Comparison of Processors}

Table \ref{tab:processor_comparison} shows the results of running GAGA on ten instances of a random graph with $n=30$ and $p=0.75$ across the two processors. We find that Processor-B is 1.3 times faster, on average, than Processor-G. All reported times in the main body of the paper are unscaled, with GAGA on Processor-G and Branch-And-Bound on Processor-B.

\begin{table}[ht]
\centering
%\begin{adjustwidth}{-1cm}{}
\renewcommand{\arraystretch}{1.1}
    \begin{tabular}{|c|c|c|c|c|c|c|c|c|}
    \hline
        \multicolumn{3}{|c|}{Instance} & \multicolumn{2}{|c|}{GAGA-SA-10000-Processor-G} & \multicolumn{2}{|c|}{GAGA-SA-10000-Processor-B} & \multicolumn{2}{|c|}{Branch-And-Bound} \\
    \hline
        n & p & i & obj & time (s) & obj & time (s) & obj & time (s) \\ \hline

    \hline
        30 & 0.75 & 0 &  1.02~e+05  &  11937.47  &  1.03~e+05  &  8927.38  & 1.17~e+05 &  14400 \\
        30 & 0.75 & 1 &  1.22~e+05  &  10765.75  &  1.23~e+05  &  8275.47  & 1.41~e+05 &  14400 \\ 
        30 & 0.75 & 2 &  0.96~e+05  &  12556.06  &  0.96~e+05  &  9534.69  & 1.20~e+05 & 14400 \\
        30 & 0.75 & 3  &  0.99~e+05  &  12773.28  &  1.00~e+05  &  9704.60  & 1.08~e+05 & 14400 \\
        30 & 0.75 & 4 &  0.95~e+05  &  12036.30  &  0.96~e+05  &  8397.12  & 1.16~e+05 & 14400 \\
        30 & 0.75 & 5 &  1.02~e+05  &  13553.03  &  1.03~e+05  &  10758.42  & 1.16~e+05 & 14400 \\
        30 & 0.75 & 6 &  0.98~e+05  &  11316.65  &  0.98~e+05  &  8588.02  & 1.13~e+05 & 14400 \\
        30 & 0.75 & 7  &  1.05~e+05  &  11515.46  &  1.06~e+05  &  8315.06  & 1.32~e+05 & 14400 \\
        30 & 0.75 & 8 &  0.94~e+05  &  11144.66  &  0.94~e+05  &  9293.79  & 1.13~e+05 & 14400 \\
        30 & 0.75 & 9 &  0.88~e+05  &  12763.55  &  0.89~e+05  &  9958.27  & 1.10~e+05 & 14400 \\
    \hline
    \end{tabular}
    %\end{adjustwidth}
    \caption{Results comparing processors on $n=30, p=0.75$.We find that Processor-B is faster by a factor of $\sim 1.3$ than Processor-G}
    \label{tab:processor_comparison}
\end{table}

\newpage

\bibliographystyle{plain}
{\small
\bibliography{reference}
}
\end{document}